\documentclass[11pt]{article}
\usepackage{amsmath}
\usepackage{amssymb}
\usepackage{hyperref}
\hsize \textwidth
\advance \hsize by -\marginparwidth
\evensidemargin \oddsidemargin
\advance\hoffset by 5mm

\newcommand{\no}{\nonumber}
\newcommand{\be}{\begin{equation}}
\newcommand{\ee}{\end{equation}}
\newcommand{\bi}{\begin{itemize}}
\newcommand{\ei}{\end{itemize}}
\newcommand{\br}{\begin{eqnarray}}
\newcommand{\er}{\end{eqnarray}}

\newcommand{\expE}{\mathbb{E}}
\newcommand{\qed}{$\Box$}

\newcommand{\abs}[1]{\lvert #1 \rvert}
\newcommand{\norm}[1]{\lVert #1 \rVert}

\newcommand{\avg}[1]{\langle #1 \rangle}

\newcommand{\davg}[1]{\langle \! \langle #1 \rangle \! \rangle}

\newcommand{\commentout}[1]{}

\newtheorem{theo}{Theorem}[section]

\newtheorem{lem}[theo]{Lemma}

\newtheorem{rmk}[theo]{Remark}

\newcommand{\Rm}{{\mathbb R}}
\newcommand{\Zm}{{\mathbb Z}}

\newcommand{\Pm}{{\mathbb P}}

\def\Zint#1{\mathchoice
{\ZZint\displaystyle\textstyle{#1}}%
{\ZZint\textstyle\scriptstyle{#1}}%
{\ZZint\scriptstyle\scriptscriptstyle{#1}}%
{\ZZint\scriptscriptstyle\scriptscriptstyle{#1}}%
\!\int}
\def\ZZint#1#2#3{{\setbox0=\hbox{$#1{#2#3}{\int}$}
\vcenter{\hbox{$#2#3$}}\kern-.5\wd0}}

\def\dashint{\Zint-}

\begin{document}

\title{Normal approximation for the net flux through a random conductor}

\author{James Nolen\thanks{Department of Mathematics, Duke University, Durham, North Carolina, USA. (nolen@math.duke.edu). }  }

\date{June 6, 2104}
\maketitle

\begin{abstract}
We consider solutions of an elliptic partial differential equation in $\Rm^d$ with a stationary, random conductivity coefficient. The boundary condition on a square domain of width $L$ is chosen so that the solution has a macroscopic unit gradient. We then consider the average flux through the domain. It is known that in the limit $L \to \infty$, this quantity converges to a deterministic constant, almost surely.  Our main result is about normal approximation for this flux when $L$ is large: we give an estimate of the Kantorovich-Wasserstein distance between the law of this random variable and that of a normal random variable. This extends a previous result of the author \cite{N1} to a much larger class of random conductivity coefficients.  The analysis relies on elliptic regularity, on bounds for the Green's function, and on a normal approximation method developed by S. Chatterjee \cite{Chat1} based on Stein's method.
\end{abstract}

\section{Introduction}

This paper pertains to solutions of the random partial differential equation
\be
- \nabla \cdot \left(a(x) (\nabla \phi(x) + e_1)\right) + \beta \phi(x) = 0, \quad x \in D_L \subset \Rm^d, \label{corPDE}
\ee
where the coefficient $a(x) = (a_{ij}(x)) \in (L^\infty(\Rm^d))^{d \times d}$ is a stationary random matrix satisfying a uniform ellipticity condition. The parameter $\beta \geq 0$ is deterministic.  The set $D_L = [0,L)^d$ is the domain, and we require that $\phi$ satisfies periodic boundary conditions on the boundary of $D_L$.  Our main result is about the statistical behavior of the quantity
\be
\Gamma_{L,\beta} =  \frac{1}{|D_L|} \int_{D_L} (\nabla \phi + e_1) \cdot a(x)(\nabla \phi + e_1)  + \beta \phi^2\,dx \label{Gammadef}
\ee
for large $L$. Using (\ref{corPDE}) and the periodicity of $\phi$ we see that $\Gamma_{L,\beta}$ may also be written as
\be
\Gamma_{L,\beta} =  \frac{1}{|D_L|} \int_{D_L} e_1 \cdot a(x)(\nabla \phi(x) + e_1) \,dx. \label{gammafluxdef}
\ee
This is a random variable, as the coefficient $a(x)$ and the solution $\phi$ are random.

Partial differential equations like (\ref{corPDE}) arise in physical applications where the coefficient $a(x)$ may be modeled best as a random field, due to inherent uncertainty and complexity of the physical medium \cite{Torq}.  If we interpret (\ref{corPDE}) in terms of electrical conductivity, then $\phi$ is a potential, $a(x)$ is the conductivity, and the vector field $-a(x) (\nabla \phi + e_1)$ is a current density. The unit vector $e_1$ is deterministic, the gradient of the linear potential $x \cdot e_1$.  Considering (\ref{gammafluxdef}), we interpret $\Gamma_{L,\beta}$ as an average flux in the direction $e_1$ that results from a macroscopic potential gradient imposed in the direction of $e_1$.  

The equation (\ref{corPDE}) plays an important role in the homogenization theory for the random elliptic operator $u \mapsto - \nabla \cdot (a(x/\epsilon) \nabla u)$ in the limit $\epsilon \to 0$ \cite{PV1, JKO}. It is well-known that the homogenized conductivity tensor $\bar a$ for that operator can be expressed in terms of functions $\phi$, called ``correctors", which solve (\ref{corPDE}) with $e_1$ being one of the $d$ standard basis vectors and which have stationary gradient. On the other hand, in a numerical computation of $\bar a$ one must approximate the true correctors by solving (\ref{corPDE}) in a bounded domain $D_L$ with suitable boundary condition. The parameter $\beta \geq 0$ is a kind of regularizing parameter that sometimes is used in approximation theory. The periodic boundary condition that we impose here is one choice that allows accurate approximation of the effective coefficient $\bar a$ in the limit $L \to \infty$ \cite{BP, Ow1, GNO}.

The results of \cite{BP, Ow1} imply that for $\beta \geq 0$ fixed, $\Gamma_{L,\beta}$ converges almost surely, as $L \to \infty$, to a deterministic constant $\bar \Gamma_\beta > 0$. For $\beta = 0$, the limit $\bar \Gamma_0$ is one of the diagonal entries of the homogenized tensor $\bar a$ described above.  For finite $L$, it is interesting to understand how $\Gamma_{L,\beta}$ and $\phi$ fluctuate around their means. Our main result is an estimate showing that for $L > \!\!> 1$, the distribution of $\Gamma_{L,\beta}$ is very close to that of a normal random variable. In \cite{N1}, we proved a similar result under strong assumptions about the random coefficient $a(x)$. In the present paper, however, we develop a more general approach which yields normal approximation for $\Gamma_{L,\beta}$ under much weaker assumptions about the law of $a(x)$.

Before we present the main result and explain its relation to other works, let us define the problem precisely and establish notation. 

\subsection*{The random coefficient $a(x)$}

For $L \in \Zm^+$, let $D_L = [0,L)^d \subset \Rm^d$ and let $L^\infty_{per}(D_L)$ denote the set of functions in $L^\infty(\Rm^d)$ which are periodic with period $L$ in each direction. That is, for all $f \in L^\infty_{per}(D_L)$, $f(x + Lk) = f(x)$ holds for all $k \in \Zm^d$ and almost every $x \in \Rm^d$. The coefficient $a(x)$ in (\ref{corPDE}) will be a random symmetric matrix with entries $a_{ij} \in L^\infty_{per}(D_L)$. Since we will be working with functions that are periodic over $D_L$, we use $dist(x,y)$ to refer to the periodized distance function:
\be
dist(x,y) = \min_{k \in \mathbb{Z}^d} d_2(x,y+k L), \quad \quad x,y \in \Rm^d, \label{distdef}
\ee
where $d_2(x,y)$ is the standard Euclidean metric in $\Rm^d$. Also, when working in the torus $D_L$, we use the notation $B_r(x)$ to refer to the ball of radius $r$ in this metric on the torus:
\br
B_r(x) = \{ y \in D_L \;|\; dist(x,y) < r \}. \label{BRdef}
\er

We suppose that the random nature of $a(x)$ comes from its dependence on a collection of $L^d$ independent random variables $Z = \{ Z_k\}_{k \in \Zm^d \cap D_L}$ taking values in a set $\mathcal{Z}$, and defined over a probability space $(\Omega,\mathcal{F},\Pm)$. Thus, $Z:\Omega \to \mathcal{Z}^{L^d}$. We often will write $a(x)$ for $a(x,Z)$, the dependence on $Z$ being understood. Let $\expE[f(Z)]$ denote expectation with respect to the probability measure $\Pm$ defining the law of $Z$. We will make three additional structural assumptions about the random matrix $a(x)$. First, we require that $a(x)$ is statistically stationary with respect to integer shifts in $x$: for every $k \in \Zm^d$ and $a(\cdot + k)$ is equal in law to $a(\cdot)$.  Second, we suppose boundedness and uniform ellipticity: there are positive constants $a^*, a_* > 0$ such that for any nonzero $\xi \in \Rm^d$
\be
a_* |\xi|^2 \leq \xi \cdot a(x)\xi \leq a^* |\xi|^2, \quad  x \in D_L \label{unifelip}
\ee
holds $\Pm$-almost surely. Third, we suppose that there is a constant $\tau > \sqrt{d} > 0$ such that for all $k \in \mathbb{Z}^d$
\br
a(x,Z) - a(x,Z')  = 0 \quad \text{if}\;\; dist(x,k) \geq \tau \label{gradalowerupper}
\er
holds whenever $Z_j = Z_j'$ for all $j \neq k$. One consequence of this last assumption is that $x \mapsto a(x,Z)$ does not depend of $Z_k$ if $dist(x,k) \geq \tau$. Moreover, $a(x,Z)$ and $a(y,Z)$ are statistically independent if $dist(x,y) \geq 2 \tau$. In other words, the dependence of $a(x,Z)$ on $Z$ is local: $a(x,Z)$ depends only on $Z_j$ for indices $j \in \mathbb{Z}^d$ that are sufficiently near $x \in \Rm^d$.

For clarity, let us highlight some simple examples for which these assumptions hold.  First, suppose that $a(x)$ is scalar and has the form of a random checkerboard
\be
a(x) = \sum_{k \in \Zm^d \cap D_L} Z_{k} \mathbb{I}_{Q_k}(x \;\text{mod}\; L) \label{astruct}
\ee
where $\{Z_k\}_{k \in \Zm^d \cap D_L}$ is a family of independent and identically distributed real-valued random variables satisfying $a_* \leq Z_k \leq a^*$ almost surely. The set $Q_k = k + [0,1)^d$ is the unit cube with a corner at $k \in \Zm^d$, and $(x \;\text{mod}\; L)$ denotes the point $(x_1 \;\text{mod}\;L, \dots,x_d \; \text{mod}\;L) \in D_L$. So, $a(x)$ is a piecewise constant function, taking random values on the cubes $Q_k$. It is also periodic over $D_L$.  In this example, $\mathcal{Z} = [a_*,a^*]$, but we need not make any further assumptions about regularity of the law of $Z_k$, as was required in \cite{N1}.

In the next example, $a(x,\omega)$ is scalar and represents pores of conductivity $a^*$ distributed randomly within a material having background conductivity $a_* > 0$. The pores are spheres having random radii, whose centers are determined by a Poisson point process with intensity $\mu > 0$. To construct such a conductivity function, let $\{X^k_j \;|\;\; k \in \mathbb{Z}^d,\; j \in \mathbb{N} \}$ be a collection of independent random variables that are each uniformly distributed on the cube $Q_0 = [0,1)^d$. Let $\{N_k\}_{k \in \mathbb{Z}^d}$ be an independent set of Poisson random variables with mean $\mu > 0$, defined on the same probability space. That is, $\Pm( N_k = n) = (n!)^{-1} e^{-\mu} \mu^n$ for $n = 0,1,2,\dots$. The random integer $N_k$ will be the number of pores with centers in the cube $Q_k = k + [0,1)^d$.  The random measure
\[
\rho_k(A) = \sum_{j=1}^{N_k} \mathbb{I}_A(k + X^k_j)
\]
on Borel sets $A \subset Q_k$ is a homogeneous Poisson point process on $Q_k$ with intensity $\mu$. Let $\{R_j^k \;|\;\; k \in \mathbb{Z}^d,\; j \in \mathbb{N} \}$ be an independent collection of identically distributed, real-valued random variables such that $\Pm(0 < R_j^k \leq R_{max}) = 1$ for some constant $R_{max}$; these are the radii of the pores. Finally, we define
\be
a(x) = a_* + (a^* - a_*) \min\left(1 \;, \;\sum_{k \in \mathbb{Z}^d \cap D_L} \sum_{j=1}^{N_k} \mathbb{I}_{B_{R_j^k}(X^k_j)}(x - k) \right). \label{poissoncoeff}
\ee
Thus, $a(x) = a^*$ if and only if 
\[
dist(x, (k + X_j^k)) < R_j^k,\quad \text{for some} \;\;j \leq N_k \;\;\text{and} \;\;k \in \mathbb{Z}^d \cap D_L.
\]
Otherwise, $a(x) = a_*$. In this case, the random variables $\{Z_k\}_{k \in \mathbb{Z}^d \cap D_L}$ are the collections of (shifted) pore centers and radii: $Z_k = \{ (X^k_j, R^k_j) \; |\; 0 \leq j \leq N_k \}$, and we may take the set $\mathcal{Z}$ to be the set of all finite sequences $\left( (x_1,r_1),\dots, (x_n,r_n) \right)$ where $x_i \in Q_0$ and $r_i \in (0,R_{max})$.  Recalling (\ref{BRdef}), we see that $a \in L^\infty_{per}(D_L)$ almost surely, and the stationarity property holds.   The condition $\Pm(R_j \leq R_{max}) = 1$ guarantees that (\ref{gradalowerupper}) holds with $\tau = R_{max} +\sqrt{d}$.  There are many variations of this construction which fit into the framework described above, such as random rods having random orientation and length, as in the experiments described in \cite{BRCCLW}.

\subsection*{The energy functional}

Let $H^1_{per}(D_L)$ denote the set of $L$-periodic functions in $H^1_{loc}(\Rm^d)$. That is, $\phi \in H^1_{per}(D_L)$ if $\phi \in H^1_{loc}(\Rm^d)$ and $\phi(x + Lk) = \phi(x)$ a.e. $\Rm^d$ for every $k \in \Zm^d$. If $a_{ij}(x) \in (L^\infty(D_L))^{d \times d}$ and satisfies (\ref{unifelip}), then there exists a weak solution $\phi \in H^1_{per}(D_L)$ to (\ref{corPDE}):
\be
\int_{D_L} \nabla v \cdot a(x) (\nabla \phi + e_1)  + \beta \phi v \,dx = 0, \;\;\;\forall \;v \in H^1_{per}(D_L). \label{weakform1}
\ee
For $\beta > 0$, the solution is unique. For $\beta = 0$, the solution is not unique, but any two solutions in $H^1_{per}(D_L)$ must differ by a constant. So, under the normalization condition
\be
\int_{D_L} \phi(x) \,dx = 0, \label{phinorm}
\ee
and for fixed $L$, the solution is unique in $H^1_{per}(D_L)$ for all $\beta \geq 0$. With $a(x) = a(x,Z)$ satisfying the conditions above, this unique solution $\phi(x) = \phi(x,a,L,\beta)$ depends on the parameters $L$ and $\beta$, on $x \in D_L$, and on the random variables $Z = \{Z_j\}_{j \in D_L \cap \Zm^d}$ which determine $a$. The uniqueness of the solution and the stationarity of $a$ implies that $\phi(x)$ is statistically stationary with respect to integer shifts: the law of $\phi(x)$ is the same as that of $\phi(x + k)$ for any $k \in \mathbb{Z}^d$.

Having defined both $a(x)$ and $\phi(x)$, we now define the random variable $\Gamma_{L,\beta}$ by (\ref{Gammadef}), which is equivalent to (\ref{gammafluxdef}). This also is a function of the $L^d$ random variables $\{Z_j\}_{j \in D_L \cap \Zm^d}$.  We will use $\Phi_j$ and $\Phi_j'$ to refer to the integrals
\be
\Phi_j = \left( \int_{Q_j} |\nabla \phi(x) + e_1|^2 \,dx \right)^{1/2}, \quad \quad \hat \Phi_j = \left( \int_{B_\tau(j)} |\nabla \phi(x) + e_1|^2 \,dx \right)^{1/2} \label{PhiPhiM}
\ee
which appear frequently in the analysis. Recall that $B_\tau(j) \supset Q_j$, so $\hat \Phi_j \geq \Phi_j$.

\commentout{
It is known that $\Gamma_{L,\beta}$ has a variational representation:
\[
\Gamma_{L,\beta} = \min_{v \in H^1_{per}(D_L)} \frac{1}{|D_L|} \int_{D_L} a(x) \abs{\nabla v + e_1}^2 + \beta v^2 \,dx.  
\]
The Euler-Lagrange equation for this variational problem is (\ref{corPDE}), and $\phi$ is the unique minimizer (unique up to addition of a constant if $\beta = 0$). 
}

\subsection*{Main result}

Our main result is the following theorem. Suppose $W$ and $Y$ are two real-valued random variables and that $\mu_W$ and $\mu_Y$ denote the laws on $\Rm$ of $W$ and $Y$, respectively.  The Kantorovich-Wasserstein distance between $\mu_W$ and $\mu_Y$ is
\br
d_{\mathcal{W}}(W,Y) & = & \sup \left\{ \; \left|\expE h(W) - \expE h(Y)  \right|\;\;|\;\; \norm{h}_{Lip} \leq 1 \right\} \no \\
& = & \sup \left\{ \; \left|\int_{\Rm} h(w) \,d\mu_W(w)  - \int_{\Rm} h(y) \,d\mu_Y(z) \right|\;\;|\;\; \norm{h}_{Lip} \leq 1 \right\}. \no
\er

\begin{theo} \label{theo:cltmain}
Let $d \geq 2$. Let $m_{L,\beta} = \expE[\Gamma_{L,\beta}]$ and $\sigma^2_{L,\beta} = Var(\Gamma_{L,\beta})$.  Let $Y$ denote a standard normal random variable, $N(0,1)$. There is a constant $C>0$ (depending only on $d$, $a_*$, and $a^*$) and a constant $q > 2$ such that
\be
d_{\mathcal{W}}\left(\frac{\Gamma_{L,\beta} - m_{L,\beta}}{\sigma_{L,\beta}},Y\right) \leq C \frac{L^{-2d}}{\sigma_{L,\beta}^3} \expE[ \Phi_0^6] + C \frac{L^{-3d/2} \log(L)}{\sigma_{L,\beta}^2} \expE[ \Phi_0^{8q}]^{\frac{1}{2q}} \label{dWZsig1}
\ee
holds for all $L > 2$ and $\beta \geq 0$.
\end{theo}

In \cite{N1}, we obtained a similar result under more restrictive structural assumptions about the law of the coefficient $a$. Specifically, the approach in \cite{N1} required that the law of $a(x)$ be obtained by a sufficiently smooth mapping of normally distributed random variables. Those assumptions excluded cases like (\ref{poissoncoeff}) where the law of $a(x)$ may have no absolutely continuous part (with respect to Lebesgue measure on $[a_*,a^*]$); the assumptions on $a(x,Z)$ in the present setting are significantly less restrictive.  Regularity of the law of $a(x)$ in \cite{N1} made it possible to differentiate $\Gamma_{L,\beta}$ with respect to the $Z_k$ and to apply a ``second order Poincar\'e inequality" developed by Chatterjee in \cite{Ch3}.  In the present setting, the more general assumptions on the law of $a(x)$ do not allow us to apply the same approach. Consequently, the proof of Theorem \ref{theo:cltmain} is based on a more general normal approximation technique from \cite{Chat1}, which is suitable for fully discrete distributions.

The variance $\sigma_{L,\beta}^2$ and the moments of the random variable $\Phi_0$ which appear in (\ref{dWZsig1}) depend on both $L$ and $\beta$.  If the moments of $\Phi$ are bounded by a constant, independent of $L$ and $\beta$, and if the variance is bounded from below by $\sigma^2_{L,\beta} \geq C L^{-d}$, then the bound (\ref{dWZsig1}) becomes
\[
d_{\mathcal{W}}\left(\frac{\Gamma_{L,\beta} - m_{L,\beta}}{\sigma_{L,\beta}},Y\right) \leq C L^{-d/2} \log L.
\]
For all dimensions $d \geq 1$, if $\beta \geq \beta_0 >  0$ is bounded away from zero independently of $L$, then all moments $\expE[\Phi_0^q]$ are bounded independently of $L > 1$ (for example, see \cite{N1}). If $\beta = 0$ or if $\beta > 0$ is allowed to vanish as $L \to \infty$, estimating the moments $\expE[\Phi_0^q]$ is a delicate issue. Elliptic regularity helps a bit.  Meyers' estimate \cite{Meyers} implies that $\nabla \phi \in L^{p^*}$ for some $p^* > 2$. If the ratio $\frac{a^*}{a_*}$ is sufficiently close to $1$, this $p^*$ may be arbitrarily large. As a result, a uniform (in $L$ and $\beta \geq 0$) bound on $\expE[ \Phi_0^q]$ follows from this regularity estimate if $\frac{a^*}{a_*} \approx 1$ (see Lemma 4.3 of \cite{N1}, for example).  This observation goes back to the work of Naddaf and Spencer \cite{NS}. On the other hand, without the assumption $\frac{a^*}{a_*} \approx 1$, the regularity only goes so far. To estimate $\expE[\Phi_0^q]$ in this situation one can use the arguments developed recently by Gloria and Otto in \cite{GO}. In that work, the authors derive variance bounds for a discrete functional similar to $\Gamma_{L,\beta}$, involving an infinite network of random resistors on the bonds of the integer lattice $\mathbb{Z}^d$. The PDE (\ref{corPDE}) is replaced by a discrete difference equation on all of $\mathbb{Z}^d$, without the periodicity assumption. The stationary potential field $\phi(x)$ is defined at points $x \in \mathbb{Z}^d$; the gradient and divergence have interpretations as difference operators. A key point in their analysis is the following bound on moments of the discrete corrector $\phi$:
\be
\expE[ \,|\phi(0)|^q \,] \leq \left\{ \begin{array}{cc} C_q , \quad \text{if}\;\;d \geq 3 \\
C_q |\log(\beta)|^{\gamma_q}, \quad \text{if}\;\;d =2. \end{array} \right. \label{GOmoment}
\ee
The constants $C_q, \gamma_q > 0$ are independent of $L > 1$ and $\beta > 0$. The analysis of \cite{GO} can be extended to the present setting (spatial continuum, with periodicity on $D_L$) to estimate moments of both $\int_{Q_0} \phi(x) \,dx$ and $\Phi_0$ (see \cite{N1} for some discussion of this). The argument shows that moments of $\int_{Q_0} \phi(x) \,dx$ satisfy the same bound as (\ref{GOmoment}), which diverges as $\beta \to 0$ if $d = 2$. On the other hand, $\Phi_0$ involves the gradient $\nabla \phi$, and it can be shown that for all $d \geq 2$, all moments $\expE[\Phi_0^q]$ are bounded independently of $L > 1$ and $\beta \geq 0$ \cite{GO3}.  In the discrete setting, the uniform control (in $L$ and $\beta$) of $\nabla \phi$ for all $d \geq 2$ was observed already by Gloria, Otto, Neukamm \cite{GNO} (see Proposition 1 therein).

In view of $\sigma_{L,\beta}^2$ appearing in (\ref{dWZsig1}), let us note that in many cases it is expected that the variance of $\Gamma_{L,\beta}$ is bounded below by $\sigma^2_{L,\beta} \geq C L^{-d}$.  Indeed, in \cite{N1} we proved that this is the case for the random checkerboard model (\ref{astruct}). This bound is closely related to earlier work of Wehr \cite{Wehr} in the discrete setting. In a forthcoming work \cite{NolenOtto1}, we will give a more general sufficient condition under which $\sigma^2_{L,\beta} \geq C L^{-d}$ holds for the continuum setting; in particular, this lower bound holds for the coefficient (\ref{poissoncoeff}) constructed from Poisson scatter.  It is not known what is the most general class of stationary random fields $a(x, \omega)$ for which the lower bound $\text{Var}(\Gamma_{L,\beta}) \geq C L^{-d}$ holds. It is conceivable that there are random fields $a(x,Z)$ satisfying both (\ref{unifelip}) and (\ref{gradalowerupper}) such that $\text{Var}(\Gamma_{L,\beta}) = o(L^{-d})$ as $L \to \infty$, due to some short-range correlation in the variables $Z_k$; the same phenomenon is possible even for simple averages of $L^d$ identically distributed random variables when the variables may be dependent on one-another.  For example, suppose $d=1$, and consider the simple case of a sequence of $L$ resistors wired together, in series, having conductivity $a_1, a_2,\dots, a_L$.  The effective conductivity of the series is just the harmonic mean $\Gamma_L =  (L^{-1} \sum_{k=1}^L \frac{1}{a_k})^{-1}$. Suppose that $\{Z_k\}$ is a sequence of independent, Bernoulli-$p$ random variables. Thus, $\Pm(Z_k = 1) = p$, $\Pm(Z_k = 0) = 1-p$. Suppose that $a_k = (2 + Z_{k+1} - Z_k)^{-1}$. These $a_k$ are dependent, but $a_k$ and $a_j$ are independent if $|k-j| > 1$. However, by definition of $a_k$, $(\Gamma_L)^{-1}$ is a telescoping sum, and the variance of the effective conductivity satisfies $\text{Var}\left( \Gamma_L \right) = O(L^{-2}) = o(L^{-d})$. Moreover, the distribution of $\Gamma_L$ (after normalization) in this simple example is not asymptotically Gaussian.

In addition to the works we have mentioned already, the two works most closely related to Theorem \ref{theo:cltmain} are those of Biskup, Salvi, and Wolff \cite{BSW} and Rossignol \cite{Ross} regarding discrete resistor network models. By making use of the martingale central limit theorem, Biskup, Salvi, and Wolff \cite{BSW} have proved a central limit theorem for a discrete quantity similar to $\Gamma_{L,\beta}$ when $\phi$ satisfies linear Dirichlet boundary conditions on a square box, in the regime of small ellipticity contrast (i.e. $|\frac{a^*}{a_*} - 1|$ is sufficiently small). Using different techniques, including generalized Walsh decomposition and concentration bounds, Rossignol \cite{Ross} has proved a variance bound and a central limit theorem for effective resistance of a resistor network on the discrete torus. We refer to the recent review paper \cite{Biskup} for many other references on the random conductance model. Also in the discrete setting, Mourrat and Otto \cite{MO14} have studied the correlation structure of the corrector itself. Delmotte and Deuschel \cite{DD} and, more recently, Marahrens and Otto \cite{MO1} derived some annealed estimates of the mixed second derivatives $\nabla_x \nabla_y G(x,y)$ of the Green function for the discrete random elliptic operator; as we mention just after Lemma \ref{lem:wksum}, there is a step in our proof which involves bounding a similar quantity.

Other works related to Theorem \ref{theo:cltmain} include those of Naddaf and Spencer \cite{NS}, Conlon and Naddaf \cite{CD2}, and Boivin \cite{BV} in the discrete case and Yurinskii \cite{VVY} in the continuum setting; they also derive upper bounds on the variance of quantities similar to $\tilde \Gamma_{L,\beta}$ and $\Gamma_{L,\beta}$. Komorowski and Ryzhik \cite{KR} have proved some related moment bounds on $\phi$ in the discrete case when $d = 1$.  If $\beta = 0$ and the dimension is $d=1$, then equation (\ref{corPDE}) can be integrated, with the solution $\phi$ written in terms of integrals of $1/a(x)$. In that case it is known that the solution itself may satisfy a central limit theorem after suitable renormalization; see Borgeat and Piatnitski \cite{BP2} Bal, Garnier, Motsch, Perrier \cite{BGMP} for precise statement of these results. In the multidimensional setting, however, those techniques do not apply.

The basis for our proof of Theorem \ref{theo:cltmain} is a normal-approximation technique of Chatterjee \cite{Chat1} (see Theorem 2.2 therein), based on Stein's method of normal approximation. This tool and related notation is explained in Section \ref{sec:normapprox}. In Section \ref{sec:determ} we give some deterministic PDE estimates (Cacciopoli's inequality and Meyers' estimate) which are used later in the analysis. Section \ref{sec:mainarg} contains the main argument in the proof of Theorem \ref{theo:cltmain}. Finally, in Section \ref{sec:Greenbound} we prove some facts about the periodic Green's function which are used in Section \ref{sec:mainarg}.

A few more comments about notation: throughout the article we will use the convention that summation over indices $j \in D_L$ means a summation over $j \in \Zm^d \cap D_L$, with $j \in \Zm^d$ being understood. For convenience we will also use brackets $\langle f \rangle = \expE[f]$ to denote expectation. We also use $C$ to denote deterministic constants that may change from line to line, but do not depend on $L$ or $\beta$. 

\vspace{0.2in}

{\bf Acknowledgment.} I am grateful to Felix Otto for very helpful discussions. This work was partially funded by grant DMS-1007572 from the US National Science Foundation.

\section{Normal approximation} \label{sec:normapprox}

In this section we summarize a general approach to normal approximation based on Stein's method, and we establish some notation that will be used throughout the paper. Suppose $W$ is a random variable with $\expE[W] = 0$ and $\expE[W^2] = 1$, and we wish to estimate 
\be
\expE h(W) - \expE h(Y) \label{hWhZ}
\ee
where $Y \sim N(0,1)$ is a standard normal random variable, $h$ is a Lipschitz continuous function on $\Rm$ and satisfying $\norm{h'}_\infty \leq 1$. Stein's method of normal approximation \cite{Stein1} is based on the following:

\begin{lem} [See \cite{Chat1}, Lemma 4.2] \label{lem:stein}
Suppose $h:\Rm \to \Rm$ is absolutely continuous with bounded derivative, and $Y \sim N(0,1)$. There exists a solution to 
\be
\psi'(x) - x \psi(x) = h(x) - \expE[h(Y)], \quad x \in \Rm \label{steineqn}
\ee
which satisfies $\norm{\psi'}_\infty \leq \sqrt{\frac{2}{\pi}} \; \norm{h'}_\infty$ and $\norm{\psi''}_{\infty} \leq 2 \norm{h'}_\infty$.
\end{lem}

Therefore, to estimate (\ref{hWhZ}) it suffices to estimate
\be
\avg{h(W) - h(Y)} = \avg{\psi'(W) - W \psi(W)} = \text{Cov}(\avg{\psi'(W)} W - \psi(W), W) \label{hcovrep}
\ee
where $\psi$ solves (\ref{steineqn}). In particular, a bound on $\text{Cov}(\avg{\psi'(W)} W - \psi(W), W)$ which is independent of $h$ satisfying $\norm{h'}_\infty \leq 1$ will imply a bound on $d_\mathcal{W}(W,Y)$.  Lemma \ref{lem:covdif} below gives us a way of estimating $\text{Cov}(\avg{\psi'(W)} W - \psi(W), W)$ when $W = f(Z_1,\dots,Z_n)$ is a function of a collection of independent random variables. Suppose $Z = (Z_1,\dots,Z_n) \in \mathcal{Z}^n$ is a random $n$-tuple in $\mathcal{Z}^n$ having components that are independent, $\mathcal{Z}$ being a given set.  Suppose $Z' = (Z_1',\dots,Z_n') \in \mathcal{Z}^n$ is an independent copy of $Z$. Define
\be
Z^{j} = (Z_1,\dots,Z_{j-1}, Z_j',Z_{j+1},\dots, Z_n). \label{Zjdef}
\ee
Similarly, for a set $A \subset \{1,\dots,n\}$, the random $n$-tuple $Z^A$ is defined by replacing $Z_\ell$ by $Z_\ell'$, for all indices $\ell \in A$. For any function $f:\mathcal{Z}^n \to \Rm$, define
\[
\Delta_j f(Z) = f(Z^j) - f(Z).
\]
This is a function of both $Z$ and $Z'$ and we sometimes write $\Delta_j f(Z,Z')$ to emphasize this point. If $j \notin A$, then define
\[
\Delta_j f(Z^A) = f(Z^{A \cup \{j\}}) - f(Z^A).
\]
Let $[n] = \{1,2,\dots,n\}$. The following identity is due to Chatterjee \cite{Chat1}: 

\begin{lem} [See \cite{Chat1} Lemma 2.3] \label{lem:covdif}
Suppose $g,f:\mathcal{Z}^n \to \Rm$ and $\avg{g(Z)^2} < \infty$, $\avg{f(Z)^2} < \infty$. Then
\be
\text{Cov}(g(Z),f(Z))= \frac{1}{2} \sum_{j=1}^n \sum_{\substack{A \subset [n]\\ j \notin A  } } K_{n,A} \avg{\Delta_j g(Z) \Delta_j f(Z^{A})}. \label{covdif}
\ee
where $K_{n,A} =  |A|! (n - |A| - 1)!/(n!)$.
\end{lem}

\vspace{0.2in}

By applying Lemma \ref{lem:covdif} with $g = f$, one can derive the well-known Efron-Stein inequality \cite{ES1, Steele1}:

\begin{lem} \label{cor:EfronStein}
Suppose $f:\mathcal{Z}^n \to \Rm$ and $\avg{f(Z)^2} < \infty$. Then
\br
\text{Var}(f(Z))  \leq \frac{1}{2} \sum_{j=1}^n \avg{|\Delta_j f(Z)|^2}. \label{efstein}
\er
\end{lem}

By applying Lemma \ref{lem:covdif} to $\text{Cov}(\avg{\psi'(W)} W - \psi(W), W)$ in (\ref{hcovrep}), one obtains the following normal approximation bound, due to Chatterjee \cite{Chat1}: 

\begin{theo} [See \cite{Chat1} Theorem 2.2] \label{theo:normalapprox}
Let $W = f(Z)$ satisfy $\avg{W} = \mu$ and $\avg{W^2} = \sigma^2$. Then 
\be
d_{\mathcal{W}} \left( \frac{W - \mu}{\sigma},Y\right) \leq  \frac{1}{2 \sigma^3} \sum_{j=1}^n \left \langle |\Delta_j f(Z)|^3 \right \rangle  + \frac{2}{\sigma^2} \text{Var}(\expE[ T(Z,Z') | Z])^{1/2}  \label{normalapproxsum}
\ee
where $Y \sim N(0,1)$ and
\[
T(Z,Z') = \frac{1}{2} \sum_{j=1}^n \sum_{\substack{A \subset [n]\\ j \notin A  } } K_{n,A}  \Delta_j f(Z) \Delta_j f(Z^{A}).
\]
\end{theo}

\commentout{
For the readers' convenience we include a proofs of Lemma \ref{lem:covdif}, Lemma \ref{cor:EfronStein}, and Theorem \ref{theo:normalapprox} in the appendix.
}

\vspace{0.2in}

Our goal will be to prove Theorem \ref{theo:cltmain} by applying Theorem \ref{theo:normalapprox} to the random variable $f(Z) = \Gamma_{L,\beta}$. The term $\text{Var}(\expE[ T(Z,Z') | Z])$ in Theorem \ref{theo:normalapprox} can be estimated by the Efron-Stein inequality (\ref{efstein}). To this end, we introduce a third $n$-tuple $Z'' = (Z_1'',Z_2'',\dots,Z_n'')$ which is an independent copy of $Z$, independent of $Z'$. Let us define
\be
Z^k = (Z_1,\dots,Z_{k-1},Z_k'',Z_{k+1},\dots). \label{Zkdef}
\ee
For any function $g(Z,Z'):\mathcal{Z}^n \times \mathcal{Z}^n \to \Rm$ we define
\be
\Delta_k g(Z,Z') = g(Z^k,Z') - g(Z,Z'). \label{Zkdef2}
\ee
In particular, $\Delta_k g(Z,Z') = 0$ if $g(Z,Z')$ does not depend on $Z_k$. We use the notation $g^k$ to denote the action of replacing $Z_k$ by $Z_k''$ in the argument of $g$:
\[
g(Z,Z')^k = g(Z^k,Z').
\]
Thus, $\Delta_k g(Z,Z') = (g(Z,Z'))^k - g(Z,Z')$. Let us emphasize that $Z^k$ will always refer to (\ref{Zkdef}) while $Z^j$ refers to (\ref{Zjdef}). The $n$-tuples $Z^j$ and $Z^k$ have the same law, but the $n$-tuple denoted by $Z^j$ is not equivalent to $Z^k$ even when the values of the indices $k$ and $j$ are the same.

Now, Lemma \ref{cor:EfronStein} implies
\[
\text{Var}(\expE[ T(Z,Z') | Z]) \leq \frac{1}{2} \sum_{k} \expE[ |\Delta_k h(Z)|^2] 
\]
where $h(Z) = \expE[ T(Z,Z') | Z]$ and 
\[
\Delta_k h(Z) =  \expE[ T(Z^k,Z') \;|\; Z^k] - \expE[ T(Z,Z') \;|\; Z]  = \expE_{Z'}[T(Z^k,Z') - T(Z,Z') \;|\;Z, Z^k].
\]
Hence,
\br
\text{Var}(\expE[ T(Z,Z') | Z]) & \leq & \frac{1}{2} \expE  \left[ \sum_{k} \expE_{Z'}[T(Z^k,Z') - T(Z,Z') \;|\;Z, Z^k]^2 \right] \no \\
& \leq & \frac{1}{2} \sum_{k}  \expE  \left[|T(Z^k,Z') - T(Z,Z') |^2 \right].
\er
Recalling that
\[
T(Z,Z') = \frac{1}{2} \sum_{j=1}^n \sum_{\substack{A \subset [n] \\ j \notin A}} K_{n,A} \Delta_j f(Z) \Delta_j f(Z^A),
\]
we conclude that
\br
\text{Var}(\expE[ T(Z,Z') | Z]) & \leq & \frac{1}{8} \sum_k \expE  \left[  \left |  \sum_{j=1}^n \sum_{\substack{A \subset [n] \\j \notin A}} K_{n,A}\left( \Delta_j f(Z) \Delta_j f(Z^A) -  (\Delta_j f(Z))^k (\Delta_j f(Z^A))^k \right)  \right|^2 \right] \no \\
& \leq & \sum_k \expE  \left[  \left |  \sum_{j=1}^n \sum_{\substack{A \subset [n] \\j \notin A}} K_{n,A} (\Delta_j f(Z) - (\Delta_j f(Z))^k) \Delta_j f(Z^A) \right|^2 \right] \no \\
& & +  \sum_k \expE  \left[  \left |  \sum_{j=1}^n \sum_{\substack{A \subset [n] \\j \notin A}} K_{n,A} (\Delta_j f(Z))^k( \Delta_j f(Z^A) - (\Delta_j f(Z^A))^k)  \right|^2 \right]. \label{bigvarbound}
\er
Let us clarify the notation here. In the case $k \neq j$, we have
\[
(\Delta_j f(Z))^k = f(Z_1,\dots,Z_j',\dots, Z_k'',\dots,Z_n) - f(Z_1,\dots,Z_j,\dots, Z_k'',\dots,Z_n) = \Delta_j f(Z^k,Z').
\]
If $k = j$, then we have $(\Delta_j f(Z))^k = \Delta_j f(Z^k,Z') = f(Z^j) - f(Z^k)$. Nevertheless, for all $j$ and $k$ we have
\[
\Delta_j f(Z) - (\Delta_j f(Z))^k = - \Delta_k \Delta_j f(Z) = - \Delta_k (\Delta_j f(Z,Z')).
\]
So, the first sum on the right side of (\ref{bigvarbound}) is
\br
\sum_k \expE  \left[  \left |  \sum_{j=1}^n  (\Delta_k (\Delta_j f(Z))) \overline{\Delta_j f(Z^A)} \right|^2 \right],  \label{firstsum}
\er
and the second sum is
\br
\sum_k \expE  \left[  \left |  \sum_{j=1}^n  (\Delta_j f(Z))^k \Delta_k \overline{ \Delta_j f(Z^A)} \right|^2 \right],  \label{secondsum}
\er
where we have used the notation $\overline{\Delta_j f(Z^A)}$ to indicate averaging with respect to the set $A$. Specifically, if $S_{n,j}$ denotes the collection of all subsets $A \subset \{1,\dots,n\}$ which do not contain the index $j$, and $H_A:S_{n,j} \to \Rm$, then 
\be
\overline{H_A} = \sum_{\substack{A \subset [n]\\ j \notin A  } } K_{n,A} H_A = \sum_{A \in S_{n,j}} K_{n,A} H_A. \label{Aavg}
\ee
The weights $K_{n,A} \geq 0$ define a probability measure on $S_{n,j}$: $\sum_{A \in S_{n,j}} K_{n,A} = 1$. 

\commentout{
To sample from this measure, choose an integer $k \in \{0,\dots,n-1\}$ uniformly at random. Having chosen $k$, pick $k$ elements uniformly at random from $\{1,\dots,n\} \setminus \{j\}$. If $k=0$, choose the empty set. Thus, the probability of choosing a particular set $A$ of size $k$ is: 
\[
\Pm(A) = \Pm(A \;|\; |A| = k ) \Pm(|A| = k) =    {n-1 \choose k}^{-1} \frac{1}{n} = K_{n,A}.
\]
}


\section{Deterministic estimates for solutions of the elliptic equation} \label{sec:determ}

In proving Theorem \ref{theo:cltmain} we will make use of some regularity estimates -- Cacciopoli's inequality and Meyers' estimate -- that apply to solutions of elliptic PDEs. These estimates rely only on the uniform ellipticity assumption, not on the statistical structure of the coefficient $a(x)$ or on the periodicity.


\subsection*{Cacciopoli's inequality} 
if $\bar u_D$ is the average of a function $u$ over a bounded domain $D$, then the Poincar\'e inequality is $\norm{u - \bar u_D}_{L^2(D)} \leq C_D \norm{\nabla u}_{L^2(D)}$. For solutions of elliptic equations, Cacciopoli's inequality gives the reverse inequality, enabling control of $\nabla u$ by $u$ itself. The basic estimate is:

\begin{lem} \label{lem:Caccio}
Let $d \geq 1$. There is a constant $K$ such that if $R > 0$ and $u \in H^1(B_R(x_0))$ is a weak solution to $- \nabla \cdot (a \nabla u) + \beta u=  \nabla \cdot \xi$ for $x \in B_R(x_0)$, with $\xi \in (L^2(B_R))^d$, then 
\br
\int_{B_{\frac{R}{2}}(x_0)}  \abs{\nabla u}^2  \, dx & \leq  & K \left(  \int_{B_R(x_0)} \abs{\xi}^2 \,dx + \frac{1}{R^2} \int_{B_R(x_0)}  (u(x) - b)^2\,dx + \beta b^2 R^d\right) \label{caccio1}
\er
holds for any constant $b \in \Rm$. 
\end{lem}

Lemma \ref{lem:Caccio} and variants are a consequence of the following:

\begin{lem} \label{lem:CaccioGen}
Let $K_1 = 2/a_*$, $K_2 = (2/a_*) + 8(a^*/a_*)^2$, and $K_3 = (2/a_*) + 2/(a_*)^2$. Let $Q$ be a bounded open subset of $\Rm^d$ with smooth boundary. If $\beta \geq 0$ and $u \in H^1(Q)$ is a weak solution to $- \nabla \cdot (a \nabla u) + \beta u =  f +  \nabla \cdot \xi$ for $x \in Q$, with  $f \in L^2(Q)$ and $\xi \in (L^2(Q))^d$, then
\br
\int_{Q} \varphi^2 \abs{\nabla u}^2  \, dx  & \leq &   K_1  \int_{Q} f (u-b) \varphi^2 \, dx  - K_1  \beta \int_{Q} u(u-b) \varphi^2 \,dx  \no \\
& & + K_2  \int_{Q} \abs{\nabla \varphi}^2 (u- b)^2\,dx  +K_3  \int_{Q} \abs{\xi}^2 \varphi^2 \,dx \label{caccioGenbound}
\er
holds for any smooth function $\varphi \geq 0$ which vanishes on the boundary of $Q$, and any constant $b \in \Rm$.
\end{lem}

\vspace{0.2in}

For proofs of Lemma \ref{lem:Caccio} and Lemma \ref{lem:CaccioGen}, see \cite{N1} (also \cite{Gia}, for example). The factor $R^{-2}$ in (\ref{caccio1}) comes from choosing a test function $\varphi$ in (\ref{caccioGenbound}) with $|\nabla \varphi| \leq R^{-1}$. There is nothing special about the balls $B_R$ and $B_{2R}$ in Lemma \ref{lem:Caccio}; for other nested domains whose boundaries are separated by distance $R$, a similar bound follows directly from Lemma \ref{lem:CaccioGen}.

\subsection*{Meyers' Estimate}

We also will make use of a well-known regularity estimate of Meyers \cite{Meyers} which shows that if $u \in H^1_{loc}$ satisfies $- \nabla \cdot(a \nabla u) + \beta u = 0$, then $\nabla u \in L^p_{loc}$ for some $p > 2$. Moreover, $\nabla u$ may be bounded as follows:

\begin{lem} \label{lem:Meyers1}
There is a constant $p^* > 2$, depending on $d$ and $a^*/a_*$, such that the following holds for all $p \in [2,p^*]$: there is $C$ such that if $R > 0$ and $u \in H^1(B_{4R}(y))$ satisfies $- \nabla \cdot(a \nabla u) + \beta u = 0$ in $B_{4R}(y)$, then
\[
\left( \dashint_{B_R} |\nabla u|^p \,dx \right)^{1/p} \leq C R^{-1}  \left( \dashint_{B_{4R}} u^2 \,dx \right)^{1/2}.
\]  
\end{lem}
{\bf Proof of Lemma \ref{lem:Meyers1}:}  This is a consequence of Theorem 2 of \cite{Meyers} and Lemma \ref{lem:CaccioGen}. Since $u$ satisfies $- \nabla \cdot a \nabla u = h$ with $h = - \beta u$, we may apply Theorem 2 of Meyers' \cite{Meyers} to $u$ (with $p_1 = 2$, $r = 2$), to conclude that for $p > 2$ sufficiently small,
\br
\left( \dashint_{B_{R}} |\nabla u|^p \right)^{1/p} & \leq & C R^{-1} \left( \dashint_{B_{2R}} |u|^2 \right)^{1/2} + C R \left( \dashint_{B_{2R}} |h|^2 \,dx \right)^{1/2} \no \\
& = &   C R^{-1} \left( \dashint_{B_{2R}} |u|^2 \right)^{1/2} + C R^{-1}  \left(R^4 \beta^2 \dashint_{B_{2R}} | u|^2 \,dx \right)^{1/2}. \label{meyers1}
\er
Now we estimate the last term in (\ref{meyers1}). Let $\varphi:\Rm^d \to [0,1]$ be a smooth function supported in $B_{3R}(y)$ and satisfying $\varphi(x) = 1$ for all $x \in B_{2R}(y)$ and satisfying $|\nabla \varphi| \leq C/R$.  Applying Lemma \ref{lem:CaccioGen} with this function $\varphi$, with $b = 0$ and with $Q = B_{3R}(y)$, we conclude that
\be
\beta \int_{B_{2R}} u^2 \,dx \leq \beta \int_{B_{3R}} u^2 \varphi^2 \,dx \leq  C \int_{B_{3R}} u^2 |\nabla \varphi|^2\,dx \leq C R^{-2} \int_{B_{3R}} u^2 \,dx. \label{u2betaR1}
\ee
Now we apply Lemma \ref{lem:CaccioGen} once more, this time in $Q = B_{4R}$, using a function $\varphi:\Rm^d \to [0,1]$ supported in $B_{4R}(y)$ and satisfying $\varphi = 1$ in $B_{3R}(y)$ and $|\nabla \varphi| \leq C/R$. We conclude
\be
\beta \int_{B_{3R}} u^2 \,dx \leq \beta \int_{B_{4R}} u^2 \varphi^2 \,dx \leq  C \int_{B_{4R}} u^2 |\nabla \varphi|^2\,dx \leq C R^{-2} \int_{B_{4R}} u^2 \,dx. \label{u2betaR2}
\ee
Combining (\ref{u2betaR1}) and (\ref{u2betaR2}) we obtain
\[
R^{4} \beta^2 \dashint_{B_{2R}} u^2 \,dx \leq C \dashint_{B_{4R}} u^2 \,dx.
\]
This combined with (\ref{meyers1}) implies the result. \hfill \qed

\section{Application to the elliptic problem} \label{sec:mainarg}

In this section we prove Theorem \ref{theo:cltmain} by applying Theorem \ref{theo:normalapprox} to the random variable $f(Z) = \Gamma_{L,\beta}(Z)$ defined by (\ref{Gammadef}). In this case, the indices $j$ in Theorem \ref{theo:normalapprox} now run over the set $\mathbb{Z}^d \cap D_L$. The first step is to compute and estimate the terms  $\Delta_j \Gamma$ and $\Delta_k \Delta_j \Gamma$ which appear in the sums (\ref{firstsum}) and (\ref{secondsum}).

\subsection*{Estimating $\Delta_j \Gamma$ and $\Delta_k \Delta_j \Gamma$}

We will make use of the following chain rule and product rule for discrete differences:
\be
\Delta_j (f(Z)^2) = (\Delta_j f(Z))(f(Z^j) + f(Z)) \label{dchain}
\ee
and
\be
\Delta_j (f(Z)g(Z)) = \frac{1}{2}(\Delta_j f(Z))(g(Z^j) + g(Z)) + \frac{1}{2}(f(Z^j) + f(Z))(\Delta_j g(Z)). \label{dprodrule}
\ee
Let us introduce the notation $a^j = a(x,Z^j)$, $a^k = a(x,Z^k)$, $\phi^j = \phi(x,Z^j)$, $\phi^k = \phi(x,Z^k)$, according to (\ref{Zjdef}) and (\ref{Zkdef}).  By the structural condition (\ref{gradalowerupper}) observe that the functions
\[
\Delta_j a = a(x,Z^j) - a(x,Z), \quad \quad \text{and} \quad \quad \Delta_k a = a(x,Z^k) - a(x,Z)
\]
are supported on the sets $B_\tau(j)$ and $B_\tau(k)$ respectively. Furthermore, (\ref{gradalowerupper}) implies that
\be
\Delta_k \Delta_j a = \Delta_k  (a(x,Z^j) - a(x,Z)) \equiv 0, \quad \text{if} \;\; dist(k,j) \geq 2 \tau. \label{DeljkAvanish}
\ee

\begin{lem} \label{lem:DeljG}
There is a constant $C$ such that
\[
L^d |\Delta_j \Gamma(Z)| \leq C (\hat \Phi_j^2(Z) + \hat \Phi_j^2(Z^j))
\]
holds for all $L > 1$, $\beta \geq 0$, $j \in \mathbb{Z}^d$, where
\[
\hat \Phi_j(Z) = \left(\int_{B_\tau(j)} |\nabla \phi(x,Z) + e_1|^2 \,dx \right)^{1/2}.
\]
Moreover, for any $q > 1$, there is $C_q$ such that 
\be
L^{qd} \expE[ |\Delta_j \Gamma(Z)|^q ] \leq C_q\expE[ |\hat \Phi_0(Z)|^{2q} ] . \label{deljfmoment}
\ee
and
\be
L^{qd} \expE[|\overline{\Delta_j \Gamma(Z^A)}|^q] \leq  C_q \expE[|\hat \Phi_0(Z)|^{2q}]. \label{GammaBarBound}
\ee
hold for all $L > 1$, $\beta \geq 0$, $j \in \mathbb{Z}^d$.
\end{lem}

\vspace{0.2in}

\begin{lem} \label{lem:DeljDelkG}
There is a constant $C$, independent of $L > 1$ and $\beta \geq 0$ such that
\br
L^{d} |\Delta_k \Delta_j \Gamma(Z)| & \leq & C\left( \hat \Phi_k^2(Z) + \hat \Phi_k^2(Z^j) + \hat \Phi_k^2(Z^{k}) + \hat \Phi_k^2(Z^{jk}) \right) \no\\
& & +  C\left( \hat \Phi_j^2(Z) + \hat \Phi_j^2(Z^j) + \hat \Phi_j^2(Z^{k}) + \hat \Phi_j^2(Z^{jk}) \right) \label{DeljDelkG1}
\er
hold for all $k, j \in \mathbb{Z}^d$. Moreover,
\br
L^{2d} |\Delta_k \Delta_j \Gamma(Z)|^2 & \leq & C \left(\hat \Phi_j^2(Z) + \hat \Phi_j^2(Z^j)+ \hat \Phi_j^2(Z^{k}) + \hat \Phi_j^2(Z^{jk})\right) \no \\
&& \quad \quad \times \int_{B_\tau(j)} |\nabla \Delta_k \phi|^2 + |\nabla \Delta_k \phi^j|^2\,dx \label{DeljDelkGamma1}
\er
holds for all $j,k$ with $dist(k,j) \geq 2\tau$.
\end{lem}

\vspace{0.2in}

{\bf Proof of Lemma \ref{lem:DeljG}:}
Using (\ref{dchain}) and (\ref{dprodrule}) and the symmetry of $a$ we compute:
\br
L^d \Delta_j \Gamma(Z) & = & \frac{1}{2} \int_{D_L} (\nabla \phi^j + e_1) \cdot (\Delta_j a)( \nabla \phi^j + e_1 ) + (\nabla \phi + e_1) \cdot (\Delta_j a)( \nabla \phi + e_1 ) \,dx \label{djgaminitial} \\
& & + \int_{D_L} \frac{1}{2}   (\nabla \Delta_j \phi) \cdot (a^j + a)(\nabla \phi^j + \nabla \phi + 2 e_1) \,dx + \beta \int_{D_L} (\Delta_j \phi)(\phi^j + \phi) \,dx. \no 
\er
Due to (\ref{weakform1}), we have
\br
&& \int_{D_L} (\nabla \Delta_j \phi)\cdot  a^j (\nabla \phi^j + e_1) + \beta (\Delta_j \phi)\phi^j \,dx = 0, \no \\
&& \int_{D_L} (\nabla \Delta_j \phi) \cdot a (\nabla \phi + e_1) + \beta (\Delta_j \phi)\phi \,dx = 0. \no
\er
Using that observation we simplify (\ref{djgaminitial}) to
\br
L^d \Delta_j \Gamma(Z) & = &  \frac{1}{2} \int_{D_L} (\nabla \phi^j + e_1) \cdot (\Delta_j a)( \nabla \phi^j + e_1 ) + (\nabla \phi + e_1) \cdot (\Delta_j a)( \nabla \phi + e_1 ) \,dx  \no \\
& &  + \frac{1}{2}\int_{D_L}  \nabla \Delta_j \phi \cdot (\Delta_j a) (\nabla \phi + e_1)  \,dx \no - \frac{1}{2}\int_{D_L}  \nabla \Delta_j \phi  \cdot (\Delta_j a) (\nabla \phi^j + e_1)  \,dx \no \\
 & = & \frac{1}{2} \int_{D_L} (\nabla \phi^j + e_1) \cdot (\Delta_j a)( \nabla \phi^j + e_1 ) + (\nabla \phi + e_1) \cdot (\Delta_j a)( \nabla \phi + e_1 ) \,dx  \no \\
& &  - \frac{1}{2}\int_{D_L} \nabla \Delta_j \phi \cdot (\Delta_j a)\nabla \Delta_j \phi\,dx \no \\
& = & \int_{D_L} (\nabla \phi^j + e_1) \cdot (\Delta_j a) (\nabla \phi + e_1) \,dx. \label{Deljf}
\er
Because $\Delta_j a = a(x,Z^j) - a(x,Z)$ vanishes outside $B_\tau(j)$ (by (\ref{gradalowerupper})), we then infer that
\br
L^d |\Delta_j \Gamma(Z)| & \leq C & \int_{B_\tau(j)} |\nabla \phi^j + e_1|^2 \,dx + C \int_{B_\tau(j)} |\nabla \phi + e_1|^2 \,dx  = C (\hat \Phi_j^2(Z) + \hat \Phi_j^2(Z^j)). \no
\er 
Since $\phi$ is stationary with respect to integer shifts and because $Z$ and $Z^j$ have the same law, the random variables $\hat \Phi_j(Z)$, $\hat \Phi_j(Z^j)$, and $\hat \Phi_0(Z)$ are identically distributed. Therefore, for any $q > 1$ there is a constant $C_q$ such that
\be
L^{qd} \expE[ |\Delta_j \Gamma(Z)|^q ] \leq C_q \expE[ |\hat \Phi_j(Z)|^{2q} ] +  C_q \expE[ |\hat \Phi_j(Z^j)|^{2q} ] = 2C_q\expE[ |\hat \Phi_0(Z)|^{2q} ],
\ee
which is (\ref{deljfmoment}). 

Now we prove (\ref{GammaBarBound}). Jensen's inequality implies
\[
|\overline{\Delta_j \Gamma(Z^A)}|^q = | \sum_{\substack{A \subset [n]\\ j \notin A }} K_{n,A}  \Delta_j \Gamma(Z^A) |^q \leq \sum_{\substack{A \subset [n]\\ j \notin A }} K_{n,A}  |\Delta_j \Gamma(Z^A)|^q .
\]
Therefore from (\ref{deljfmoment}) we obtain
\[
L^{qd} \expE[|\overline{\Delta_j \Gamma(Z^A)}|^q] \leq L^{qd} \sum_{\substack{A \subset [n]\\ j \notin A }} K_{n,A} \expE[|\Delta_j \Gamma(Z^A)|^q] \leq C \expE[|\hat \Phi_0(Z)|^{2q}].
\]
\hfill \qed

\vspace{0.2in}

{\bf Proof of Lemma \ref{lem:DeljDelkG}:} Starting from (\ref{Deljf}) and using (\ref{dchain}) and (\ref{dprodrule}) we compute
\br
L^d \Delta_k \Delta_j \Gamma(Z) & = & \frac{1}{2} \int_{D_L} (\nabla \phi^{jk} + e_1) \cdot  (\Delta_k \Delta_j a)(\nabla \phi^k + e_1) + (\nabla \phi^j + e_1)\cdot (\Delta_k \Delta_j a) (\nabla \phi + e_1) \,dx \no \\ 
& & + \frac{1}{4} \int_{D_L}  \Delta_k \nabla \phi^j \cdot ((\Delta_j a)^k + (\Delta_j a))  (\nabla \phi^{k} +e_1 +  \nabla \phi + e_1) \no \\
& & + \frac{1}{4} \int_{D_L}  \Delta_k \nabla \phi \cdot ((\Delta_j a)^k + (\Delta_j a)) (\nabla \phi^{jk} + e_1 + \nabla \phi^j + e_1). \no 
\er
The matrices $\Delta_j a$ and $\Delta_k a$ are zero outside $B_\tau(j)$ and $B_\tau(k)$, respectively. Also, $\Delta_k \nabla \phi = \Delta_k (\nabla \phi + e_1) = (\nabla \phi^k + e_1) - (\nabla \phi + e_1)$. Therefore, by the Cauchy-Schwarz inequality we obtain
\br
L^{d} |\Delta_k \Delta_j \Gamma(Z)| & \leq & C\left( \hat \Phi_k^2(Z) + \hat \Phi_k^2(Z^j) + \hat \Phi_k^2(Z^{k}) + \hat \Phi_k^2(Z^{jk}) \right) \no\\
& & +  C\left( \hat \Phi_j^2(Z) + \hat \Phi_j^2(Z^j) + \hat \Phi_j^2(Z^{k}) + \hat \Phi_j^2(Z^{jk}) \right). \label{jkequal} 
\er
for all $j,k \in \mathbb{Z}^d$. This is (\ref{DeljDelkG1}). 

If $dist(k,j) \geq 2\tau$ then $(\Delta_k \Delta_j a) \equiv 0$ and $(\Delta_j a)^k = \Delta_j a$, by (\ref{DeljkAvanish}).  So, in this case we have
\br
L^d \Delta_k \Delta_j \Gamma(Z) & = &  \frac{1}{2} \int_{D_L} (\Delta_k \nabla \phi^j) \cdot  (\Delta_j a) (\nabla \phi^{k} +e_1 +  \nabla \phi + e_1) \no \\
& & + \frac{1}{2} \int_{D_L} (\Delta_k \nabla \phi) \cdot ( \Delta_j a) (\nabla \phi^{jk} + e_1 + \nabla \phi^j + e_1). \no 
\er
Applying Cauchy-Schwarz to this, using the fact that $\Delta_j a$ is supported on $B_\tau(j)$, we obtain
\br
L^{2d} |\Delta_k \Delta_j \Gamma(Z)|^{2}  & \leq & C  \int_{B_\tau(j)} |\nabla (\Delta_k  \phi^j)|^2 \,dx \left( \int_{B_\tau(j)} |\nabla \phi^k + e_1|^2 \,dx +  \int_{B_\tau(j)} |\nabla \phi + e_1|^2 \,dx \right) \no \\
& & + C  \int_{B_\tau(j)} |\nabla (\Delta_k  \phi)|^2 \,dx \left( \int_{B_\tau(j)} |\nabla \phi^{jk} + e_1|^2 \,dx +  \int_{B_\tau(j)} |\nabla \phi^j + e_1|^2 \,dx \right) \no \\\label{Delkjf}
\er
if $dist(k,j) \geq 2\tau$.  The implies (\ref{DeljDelkGamma1}).
\hfill \qed

\subsection*{Relation to the periodic Green's function}

The function $w_k = \Delta_k \phi = \phi^k - \phi \in H^1_{per}(D_L)$ which appears in Lemma \ref{lem:DeljDelkG} satisfies the equation
\be
- \nabla \cdot (a \nabla w_k) + \beta w_k = \nabla \cdot (\Delta_k a) (\nabla \phi^k + e_1), \label{DeltakEqn}
\ee
and the distribution on the right side of (\ref{DeltakEqn}) is supported on $\overline{B_\tau(k)}$. Choosing $w_k$ itself as a test function for (\ref{DeltakEqn}), we obtain the bound
\be
\int_{D_L} |\nabla w_k |^2 \,dx \leq \left(\frac{a^*}{a_*}\right)^2 \int_{B_\tau(k)} |\nabla \phi^k + e_1|^2 \,dx = \left(\frac{a^*}{a_*}\right)^2 \hat \Phi_k^2(Z^k). \label{wkalways}
\ee
Later it will be convenient to normalize the function $w_k = \Delta_k \phi$ by defining
\be
\tilde w_k =  \hat \Phi_k(Z^k)^{-1} w_k =  \left( \int_{B_\tau(k)} |\nabla \phi^k + e_1|^2 \,dx \right)^{-1/2} w_k. \label{wkdef}
\ee
The following estimate relates $w_k$ to the periodic Green's function, and it will enable us to control the decay of $|\nabla w_k|^2$ away from $B_{\tau}(k)$ (using Cacciopoli's inequality). This connection between the Green's function and quantities analogous to $\Delta_k \phi$ has been used in other works, as well (e.g. \cite{NS, GO, G1}).

\begin{lem} \label{lem:wGrep} Let $d \geq 1$, and let $G = G(x,y,Z)$ be the periodic Green's function associated with the coefficient $a(x,Z)$:
\[
- \nabla_x \cdot (a(x,Z) \nabla_x G) + \beta G = \delta_y(x) - |D_L|^{-1},
\]
normalized by $\int_{D_L} G(x,y) \,dx = 0$ in the case $\beta = 0$. There is a constant $C$ (depending only on $d,a_*,a^*$) such that for any $L > 2$, any $k \in D_L \cap \Zm^d$, and any open set $A \subset D_L$ with $dist(A,B_\tau(k)) > 0$, we have
\be
\int_{A} (\Delta_k \phi)^2 \,dy \leq  C \hat \Phi_k^2(Z^k) \int_{y \in A} \int_{x \in B_\tau(k)} |\nabla_x G(x,y)|^2 \,dx \,dy. \label{wGint}
\ee
with probability one.
\end{lem}

{\bf Proof of Lemma \ref{lem:wGrep}:} Let us define $\xi_k = (\Delta_k a) (\nabla \phi^k + e_1)$ which is supported in $B_\tau(k)$. Let $v \in H^1_{per}(D_L)$ satisfy
\[
- \nabla \cdot( a \nabla v) + \beta v = \Delta_k \phi \mathbb{I}_{A}(x) - \frac{1}{D_L}\int_{A} \Delta_k \phi(x) \,dx.
\]
By using (\ref{DeltakEqn}) and the fact that $\int_{D_L} \Delta_k \phi(x) \,dx = 0$, we have
\br
\int_{A}(\Delta_k \phi(x))^2\,dx & = &\int_{D_L} (\mathbb{I}_{A}(x) \Delta_k \phi(x))\Delta_k \phi(x) \,dx  \no \\
& = & \int_{D_L} \nabla v \cdot a(x) \nabla (\Delta_k \phi)  + \beta v (\Delta_k \phi) \,dx \no \\
& = & - \int_{D_L} \xi_k(x) \cdot \nabla v(x) \,dx \leq \left( \int_{B_\tau(k)} |\xi_k|^2 \right)^{1/2}\left( \int_{B_\tau(k)} |\nabla v|^2 \right)^{1/2}. \label{wxiv}
\er
On the other hand,
\[
v(x) = \int_{A} G(x,y) \Delta_k \phi(y) \,dy, \quad \quad  \quad \quad \nabla v(x) =  \int_{A} \nabla_x G(x,y)\Delta_k \phi(y) \,dy
\]
hold for almost every $x$ outside $A$. Therefore, by Cauchy-Schwarz we have
\[
|\nabla v(x)|^2 \leq   \int_{A} |\nabla_x G(x,y)|^{2} \,dy  \int_{A} (\Delta_k \phi(y))^2 \,dy . 
\]
for almost every $x$ in $B_\tau(j)$. Also, $\int_{B_\tau(k)} |\xi_k|^2 \,dx \leq C_2^2 \hat \Phi_k$, by (\ref{gradalowerupper}). Combining this with (\ref{wxiv}) we obtain (\ref{wGint}). \qed

\vspace{0.2in}

\vspace{0.2in}

In view of Lemma \ref{lem:DeljDelkG} and Lemma \ref{lem:wGrep}, we see that estimates of the Green's function will play an important role in estimating $\Delta_k \Delta_j \Gamma$.  We will make use of the following bounds, proved later in Section \ref{sec:Greenbound}. The first is a bound on the decay of $G(x,y)$ which is uniform with respect to the probability measure $\Pm$. The second, is a version of Lemma 2.9 in \cite{GO}, and it is also uniform with respect to the probability measure $\Pm$. Recall the definitions (\ref{distdef}) and (\ref{BRdef}) of $dist(x,y)$ and $B_r(x)$.

\begin{lem} \label{lem:udecay} 
Let $d \geq 3$. There is a constant $C > 0$, depending only on $d$, $a^*$, and $a_*$, such that
\[
|G(x,y)| \leq C \; dist(x,y)^{2-d}
\]
holds for all $x,y \in D_L$ with $x \neq y$, all $L \geq 1$ and $\beta \geq 0$. 
\end{lem}

\begin{lem} \label{lem:ugrad2d} 
Let $d =2$. There is a constant $C > 0$, depending only on $d$, $a^*$, and $a_*$, such that for all $R > 0$, $\beta \geq 0$, $L \geq 1$,
\[
\int_{B_R(x_0)} |\nabla_x G(x,y)|^2 \,dx \leq C
\]
holds for all $x_0 \in D_L$ and $y \in D_L \setminus \overline{B_{2R}(x_0)}$. 
\end{lem}

\vspace{0.2in}

\subsection*{Proof of Theorem \ref{theo:cltmain}.}

Because of the stationarity assumption, moments of $\hat \Phi_0$ are controlled by the same moments of $\Phi_0$: for any $q \geq 1$ there is a constant $C_q$ such that
\be
\expE[ |\hat \Phi_0|^q] \leq C_q \expE[ | \Phi_0|^q] \label{momentequiv}
\ee
for all $L \geq 1$ and $\beta \geq 0$. This is proved in Lemma 4.2 of \cite{N1}, for example. Therefore, according to Lemma \ref{lem:DeljG}, we can bound the first term on the right side of (\ref{normalapproxsum}) as
\be
\frac{1}{2\sigma^3} \sum_{j \in D_L} \langle |\Delta_j \Gamma(Z)|^3 \rangle \leq C \frac{L^{d}}{ \sigma^3} L^{-3d}\expE[ |\hat \Phi_0(Z)|^{6} ] \leq C \frac{L^{d}}{ \sigma^3} L^{-3d}\expE[ |\Phi_0(Z)|^{6} ]. \label{sig3sum}
\ee

As shown already, the term $\text{Var}(\expE[T(Z,Z')|Z])$ in (\ref{normalapproxsum}) is controlled by the sum of (\ref{firstsum}) and (\ref{secondsum}). We now focus on estimating (\ref{firstsum}). By Minkowsi's inequality we have
\br
\Big \langle  \left |  \sum_{j\in D_L}  (\Delta_k \Delta_j \Gamma(Z)) \overline{\Delta_j \Gamma(Z^A)} \right|^2 \Big \rangle \leq \left( \sum_{j \in D_L}  \langle |\Delta_k \Delta_j \Gamma(Z) \overline{\Delta_j \Gamma(Z^A)}|^2  \rangle^{1/2} \right)^2. \label{minkow1}
\er
It will be convenient to split up this sum over domains resembling dyadic annuli centered around the cube $Q_k = k + [0,1)^d$. Let $N$ denote the smallest integer such that $2^{N} \tau \geq L/4$. Hence, $N = O(\log L)$ and $L/4 \leq 2^{N} \tau \leq L/2$. Then, let $A_0^k$ denote the union of cubes that are close to $Q_k$:
\[
A_0^k =  \{ x \in D_L \;|\;\; x \in Q_j, \;\; 0 \leq dist(Q_j,Q_k) \leq 2 \tau \},
\]
and for $\ell = 1,2,\dots,N-1$ let $A_\ell^k$ denote the set
\[
A_\ell^k = \{ x \in D_L \;|\;\; x \in Q_j, \;\; 2^{\ell}\tau < dist(Q_j,Q_k) \leq 2^{\ell + 1}\tau \}.
\]
Again, we use $dist(Q_j,Q_k)$ to refer to distance on the torus $D_L$ (modulo $L \mathbb{Z}^d$) between sets $Q_j$ and $Q_k$. Finally, define $A_N^k$ by
\[
A_N^k = \{ x \in D_L \;|\;\; x \in Q_j, \;\; 2^{N}\tau < dist(Q_j,Q_k) \}.
\]
Each set $A_\ell^k$ is a union of cubes, and has Lebesgue measure $|A_\ell^k| = O(2^{\ell d})$. Let $A^k_+ = D_L \setminus A_0^k$. Observe that $j \in A_\ell^k$ if and only if $Q_j \subset A_\ell^k$. Similarly, $j \in A_+^k$ if and only if $Q_j \in A_\ell^k$ for some $\ell \geq 1$. Thus, 
\[
D_L = A^k_0 \cup A^k_+ = \bigcup_{\ell = 0}^N A_\ell^k.
\]
In this way, we write the sum appearing in (\ref{minkow1}) as:
\br
\sum_{j \in D_L} \langle |\Delta_k \Delta_j \Gamma(Z) \overline{\Delta_j \Gamma(Z^A)}|^2 \rangle^{1/2}  & = &   \sum_{j \in A^k_0 \cup A^k_1}  \langle |\Delta_k \Delta_j \Gamma(Z) \overline{\Delta_j \Gamma(Z^A)}|^2 \rangle^{1/2} \no \\
& &  +  \sum_{\ell = 2}^N \sum_{j \in A^k_\ell}  \langle |\Delta_k \Delta_j \Gamma(Z) \overline{\Delta_j \Gamma(Z^A)}|^2 \rangle^{1/2}. \label{minkow2}
\er
We will bound the terms in (\ref{minkow2}) using the following Lemma.
The first estimate will bound the terms with indices $j \in A_0^k \cup A^k_1$. The second estimate will be used for the other indices.

\begin{lem} \label{lem:genHolder}
For $p > 1$ there is a constant $C_p$ such that if $k,j \in D_L$ with $dist(k,j) \geq 2 \tau$ then
\br
&& L^{4d} \langle |\Delta_k \Delta_j \Gamma(Z) \overline{\Delta_j \Gamma(Z^A)}|^2\rangle  \leq   C_p \langle \hat \Phi_0^{4q} \rangle^{\frac{4}{2q}} \left\langle (\int_{B_\tau(j)} |\nabla \tilde w_k|^{2} \,dx)^p \right \rangle^{1/p} 
\er
where $q = 2p/(p-1)$. Also, there is a constant $C$ such that 
\be
L^{4d} \langle |\Delta_k \Delta_j \Gamma(Z) \overline{\Delta_j \Gamma(Z^A)}|^2 \rangle \leq  C \langle \hat \Phi_0^8\rangle \label{anykjfbound}
\ee
holds for all $j,k \in D_L$, $L > 1$, $\beta \geq 0$. 
\end{lem}
{\bf Proof:} First we prove (\ref{anykjfbound}). By Lemma \ref{lem:DeljDelkG}, we always have
\br
L^{2d} |\Delta_k \Delta_j \Gamma(Z)|^2 & \leq & C \left(\hat \Phi_j^2(Z) + \hat \Phi_j^2(Z^j)+ \hat \Phi_j^2(Z^{k}) + \hat \Phi_j^2(Z^{jk})\right) \no \\
&& \quad \quad \times \left(\hat \Phi_k^2(Z) + \hat \Phi_k^2(Z^j)+ \hat \Phi_k^2(Z^{k}) + \hat \Phi_k^2(Z^{jk})\right). \label{DeljDelkG2w2}
\er
Moreover, the terms $\hat \Phi_j^2(Z)$, $\hat \Phi_j^2(Z^j)$, $\hat \Phi_j^2(Z^{k})$, $\hat \Phi_j^2(Z^{jk})$, $\hat \Phi_k^2(Z)$, $\hat \Phi_k^2(Z^{j})$, $\hat \Phi_k^2(Z^k)$, $\hat \Phi_k^2(Z^{jk})$ are identically distributed, all having the same distribution as $\hat \Phi_0(Z)$. By Lemma \ref{lem:DeljG}, we know that 
\be
L^{2d} \langle  | \overline{\Delta_j \Gamma(Z^A)}|^{2q} \rangle^{1/q}  \leq C \langle \hat \Phi_0^{4q} \rangle^{\frac{1}{q}}. \label{GamBar2q}
\ee
Therefore, by applying the generalized H\"older inequality with $\frac{1}{4} + \frac{1}{4} + \frac{1}{2} = 1$ we obtain
\br
L^{2d} \langle |\Delta_k \Delta_j \Gamma(Z) \overline{\Delta_j \Gamma(Z^A)}|^2\rangle   &\leq &   C \left \langle \left(\hat \Phi_j^2(Z) + \hat \Phi_j^2(Z^j)+ \hat \Phi_j^2(Z^{k}) + \hat \Phi_j^2(Z^{jk})\right) \right.\no \\
&& \quad \quad \times \left. \left(\hat \Phi_k^2(Z) + \hat \Phi_k^2(Z^j)+ \hat \Phi_k^2(Z^{k}) + \hat \Phi_k^2(Z^{jk})\right) |\overline{\Delta_j \Gamma(Z^A)}|^2 \right \rangle \no \\
& \leq & C \langle \hat \Phi_0^{8} \rangle^{\frac{1}{4}} \langle \hat \Phi_0^{8} \rangle^{\frac{1}{4}}  \langle  | \overline{\Delta_j \Gamma(Z^A)}|^{4} \rangle^{1/2} \no \\
& \leq & C \langle \hat \Phi_0^{8} \rangle^{\frac{1}{4}} \langle \hat \Phi_0^{8} \rangle^{\frac{1}{4}}  L^{-2d} \langle \hat \Phi_0^{8} \rangle^{\frac{1}{2}}.  \label{genHolderbound}
\er
This proves (\ref{anykjfbound}).

If $dist(k,j) \geq 2 \tau$, Lemma \ref{lem:DeljDelkG} tells us that
\br
L^{2d} |\Delta_k \Delta_j \Gamma(Z)|^2 & \leq & C \left(\hat \Phi_j^2(Z) + \hat \Phi_j^2(Z^j)+ \hat \Phi_j^2(Z^{k}) + \hat \Phi_j^2(Z^{jk})\right) \no \\
&& \quad \quad \times \left(  \hat \Phi_k^2(Z^k) \int_{B_\tau(j)} |\nabla \tilde w_k|^2 \,dx +  \hat \Phi_k^2(Z^{jk}) \int_{B_\tau(j)} |\nabla \tilde w_k^j|^2\,dx \right), \label{DeljDelkG2w}
\er
where 
\[
\tilde w_k = \hat \Phi_k(Z^k)^{-1} \Delta_k \phi, \quad \quad \tilde w_k^j = \hat \Phi_k(Z^{jk})^{-1} \Delta_k \phi^j.
\]
Let $p > 1$, let $q = 2p/(p-1)$ so that $\frac{1}{p}  + \frac{1}{2q} + \frac{1}{2q} + \frac{1}{q}= 1$. Then by (\ref{DeljDelkG2w}) and the generalized H\"older inequality,
\br
&& L^{2d} \langle |\Delta_k \Delta_j \Gamma(Z) \overline{\Delta_j \Gamma(Z^A)}|^2\rangle   \no \\
& & \quad \quad \quad \quad \leq   C \langle \hat \Phi_0^{4q} \rangle^{\frac{1}{2q}} \langle \hat \Phi_0^{4q} \rangle^{\frac{1}{2q}} \langle (\int_{B_\tau(j)} |\nabla \tilde w_k|^{2} \,dx)^p \rangle^{1/p} \langle  | \overline{\Delta_j \Gamma(Z^A)}|^{2q} \rangle^{1/q} \no \\
& & \quad \quad \quad \quad \quad +  C \langle \hat \Phi_0^{4q} \rangle^{\frac{1}{2q}} \langle \hat \Phi_0^{4q} \rangle^{\frac{1}{2q}} \langle (\int_{B_\tau(j)} |\nabla \tilde w_k^j|^{2} \,dx)^p \rangle^{1/p} \langle  | \overline{\Delta_j \Gamma(Z^A)}|^{2q} \rangle^{1/q}.  \label{genHolderbound2}
\er
If $j \neq k$, then $\Delta_k \phi$ and $\Delta_k \phi^j$ have the same distribution (since $(Z,Z^k)$ and $(Z^j,Z^{jk})$ have the same joint distribution). Similarly, $\tilde w_k$ and $\tilde w_k^j$ must have the same distribution. Therefore,
\[
\Big \langle (\int_{B_\tau(j)} |\nabla \tilde w_k^j|^{2} \,dx)^p  \Big \rangle^{1/p} =   \Big \langle (\int_{B_\tau(j)} |\nabla \tilde w_k|^{2} \,dx)^p \Big \rangle^{1/p}.
\]
holds for all $j \neq k$. Combining this observation with (\ref{genHolderbound2}) and (\ref{GamBar2q}) we obtain
\br
&& L^{4d} \langle |\Delta_k \Delta_j \Gamma(Z) \overline{\Delta_j \Gamma(Z^A)}|^2\rangle  \leq   C \langle \hat \Phi_0^{4q} \rangle^{\frac{4}{2q}} \Big \langle (\int_{B_\tau(j)} |\nabla \tilde w_k|^{2} \,dx)^p \Big \rangle^{1/p}.
\er
This completes the proof of Lemma \ref{lem:genHolder}. \hfill \qed

\vspace{0.2in}

Now we return to (\ref{minkow2}). For the first sum on the right side of (\ref{minkow2}), over indices $j$ near $k$, we apply Lemma \ref{lem:genHolder} to obtain 
\br
L^{2d} \sum_{j \in A^k_0 \cup A^k_1} \langle |\Delta_k \Delta_j \Gamma(Z) \overline{\Delta_j \Gamma(Z^A)}|^2 \rangle^{1/2} \leq C  (|A^k_0| + |A^k_1|) \langle \hat \Phi_0^8 \rangle^{1/2}. \label{genHolderapp0}
\er
For the second sum in (\ref{minkow2}), we apply Lemma \ref{lem:genHolder} again to obtain
\br
L^{2d} \sum_{\ell=2}^N \sum_{j \in A^k_\ell} \langle |\Delta_k \Delta_j \Gamma(Z) \overline{\Delta_j \Gamma(Z^A)}|^2 \rangle^{1/2}  & \leq & C \langle \hat \Phi_0^{4q} \rangle^{\frac{1}{q}} \sum_{\ell=2}^N \sum_{j \in A^k_\ell}  \Big \langle \int_{B_\tau(j)} |\nabla \tilde w_k|^{2p} \,dx \Big \rangle^{\frac{1}{2p}}. \label{genHolderapp1}
\er
From our definition of the annuli $A^k_\ell$ and $\tau > \sqrt{d}$, we see that
\[
\bigcup_{ \substack{j \in A_\ell^k \\ \ell \geq 2}} B_\tau(j) \subset \bigcup_{ \substack{j \in A_\ell^k \\ \ell \geq 1}} Q_j.
\]
Furthermore, each ball $B_{\tau}(j)$ intersects only finitely many cubes ($O(\tau^d)$ of them). So, the last integral in (\ref{genHolderapp1}) can be replaced by an integral over $Q_j$, at the expense of a constant factor of order $O(\tau^d)$. Indeed, by Minkowski's inequality,
\br
\Big \langle \int_{B_\tau(j)} |\nabla \tilde w_k|^{2p} \,dx \Big \rangle^{\frac{1}{2p}} \leq \sum_{\substack{n \in D_L\\ |B_{\tau}(j) \cap Q_n| > 0}}  \Big \langle \int_{Q_n} |\nabla \tilde w_k|^{2p} \,dx \Big \rangle^{\frac{1}{2p}}.
\er
Therefore, (\ref{genHolderapp1}) yields
\br
L^{2d} \sum_{\ell=2}^N \sum_{j \in A^k_\ell}  \langle |\Delta_k \Delta_j \Gamma(Z) \overline{\Delta_j \Gamma(Z^A)}|^2 \rangle^{1/2} & \leq & C \langle \hat \Phi_0^{4q} \rangle^{\frac{1}{q}}\sum_{\ell=2}^N \sum_{j \in A^k_\ell}  \sum_{\substack{n \in D_L\\ |B_{\tau}(j) \cap Q_n| > 0}}  \Big \langle \int_{Q_n} |\nabla \tilde w_k|^{2p} \,dx \Big \rangle^{\frac{1}{2p}} \no \\
 & \leq & C \langle \hat \Phi_0^{4q} \rangle^{\frac{1}{q}} \sum_{\ell=2}^N \sum_{j \in A^k_\ell}   \sum_{\substack{n \in D_L\\ dist(Q_j,Q_n) < \tau}}  \Big \langle \int_{Q_n} |\nabla \tilde w_k|^{2p} \,dx \Big \rangle^{\frac{1}{2p}} \no \\
& \leq &   C \tau^d \langle \hat \Phi_0^{4q} \rangle^{\frac{1}{q}} \sum_{\ell=1}^N \sum_{j \in A^k_\ell}  \Big \langle \int_{Q_j} |\nabla \tilde w_k|^{2p} \,dx \Big \rangle^{\frac{1}{2p}} \no \\
& = &   C \tau^d \langle \hat \Phi_0^{4q} \rangle^{\frac{1}{q}} \sum_{j \in A^k_+}  \Big \langle \int_{Q_j} |\nabla \tilde w_k|^{2p} \,dx \Big \rangle^{\frac{1}{2p}}.  \label{genHolderapp}
\er

We will now show that the last sum in (\ref{genHolderapp}) is $O(\log L)$. 

\begin{lem} \label{lem:wksum}
There is are constants $C > 0$ and $p > 1$ such that
\be
 \sum_{j \in A^k_+} \Big \langle \int_{Q_j} |\nabla \tilde w_k|^{2p} \,dx \Big \rangle^{\frac{1}{2p}} \leq C \log L \label{wksumbound}
\ee
and
\be
 \sum_{k \in A^0_+}  \Big \langle \int_{Q_0} |\nabla \tilde w_k|^{2p} \,dx \Big \rangle^{\frac{1}{2p}} \leq C \log L \label{wksumbound2}
\ee
for all $\beta \geq 0$, $L \geq 2$, $k \in D_L \cap \mathbb{Z}^d$, where $\tilde w_k$ is defined by (\ref{wkdef}). 
\end{lem}

\vspace{0.2in}

Lemma \ref{lem:wGrep} gives control of  $\tilde w_k(y)$ in terms of $\nabla_x G(x=k,y)$. So, thinking heuristically, we expect that for $dist(y,k) \gg 1$, $\nabla_y \tilde w_k(y)$ should decay like the mixed second derivative $\nabla_y \nabla_x G(k,y)$ of the Green function. So, if the constant-coefficient case is any guide, we should hope that $\nabla_y \tilde w_k(y)$ decays like $O(|y - k|^{-d})$. Although we do not have uniform pointwise bounds on $\nabla_y \tilde w_k(y)$ of this sort, we still obtain (\ref{wksumbound}), which is what we would obtain if we did have the uniform bound $|\nabla_y \tilde w_k(y)| \leq C(1 + |y - k|)^{-d}$.  In the proof below, the strategy is to use Cacciopoli's inequality to control $\nabla \tilde w_k$ by $\tilde w_k$, then Lemma \ref{lem:wGrep} to control $\tilde w_k$ by $\nabla G$. Then we use stationarity and Cacciopoli's inequality again to control $\nabla G$ by $G$, for which we have uniform bounds in Lemma \ref{lem:udecay} ($d \geq 3$). Cacciopoli's inequality is applied over a large domain (the dyadic annuli) to take advantage of the $R^{-2}$ factor in Lemma \ref{lem:Caccio}.  In the context of the discrete version of this elliptic problem, a similar strategy is employed by Gloria and Otto \cite{GO} to control the decay of $\nabla_x G(x,y)$ in terms of the uniform decay of $G(x,y)$ and by Marahrens and Otto \cite{MO1} to estimate moments $\avg{|\nabla_x \nabla_y G(x,y)|^{2p}}^{1/(2p)} \leq O( (1 + |x - y|)^{-d})$ of the discrete second derivative of $G$.

\vspace{0.2in}

{\bf Proof of Lemma \ref{lem:wksum}.}  By stationarity, we have 
\[
\Big \langle \int_{Q_j} |\nabla \tilde w_k|^{2p} \,dx \Big \rangle^{\frac{1}{2p}} =  \Big \langle \int_{Q_0} |\nabla \tilde w_{k-j}|^{2p} \,dx \Big \rangle^{\frac{1}{2p}},
\]
so the bound (\ref{wksumbound2}) is equivalent to (\ref{wksumbound}). Therefore, we focus on proving (\ref{wksumbound}).

The constant $p > 1$ may be chosen so that $2p \in (0,p^*)$, where $p^* > 2$ is as in Lemma \ref{lem:Meyers1}.   We split the (\ref{wksumbound}) over the diadic annuli, and apply H\"older's inequality with $2p$ and $\frac{2p}{2p-1}$:
\br
\sum_{\ell = 1}^N \sum_{j \in A^k_\ell}  \Big \langle \int_{Q_j} |\nabla \tilde w_k|^{2p} \,dx \Big \rangle^{\frac{1}{2p}} & \leq &  \sum_{\ell = 1}^{N} \left( \sum_{j \in A_\ell^k} 1^{2p/(2p-1)} \right)^{(2p-1)/(2p)}  \left(\sum_{j \in A_\ell^k}  \Big\langle \int_{Q_j} |\nabla \tilde w_k|^{2p} \,dx  \Big \rangle \right)^{1/{2p}} \no \\
& = &   \sum_{\ell=1}^{N} |A_\ell^k|^{(2p-1)/(2p)} \left( \Big \langle \int_{A_\ell^k} |\nabla \tilde w_k|^{2p} \,dx \Big \rangle \right)^{1/(2p)}. \label{Holdersum1}
\er
For $\ell \geq 1$, let us use the notation $2A_\ell^k$ to refer to the fattened annuli:
\[
2A_\ell^k = \{ x \in D_L \;|\; x \in  Q_j,\;\; 2^{\ell-1} \tau  < dist(Q_j,Q_k) \leq 3 \cdot 2^{\ell} \tau \}, \quad \ell = 1,\dots,N-1,
\]
and
\[
2A_N^k = \{ x \in D_L \;|\; x \in  Q_j,\;\; 2^{N-1} \tau < dist(Q_j,Q_k) \}.
\]
Observe that $A^k_\ell \subset 2A^k_\ell$ and $dist(A^k_\ell, \partial( 2A^k_\ell)) \geq C 2^\ell$. Also, $dist(2 A^k_\ell, B_\tau(k)) > 0$. By Lemma \ref{lem:Meyers1} applied to $\tilde w_k$ and by Lemma \ref{lem:wGrep}, we know that
\br
\int_{A_\ell^k} |\nabla \tilde w_k|^{2p} \,dx & \leq & C (2^{\ell})^{d - p(2 + d)}  \left( \int_{2A_\ell^k} (\tilde w_k)^{2} \,dy \right)^p \no \\
& \leq & C (2^{\ell})^{d - p(2 + d)} \left(\int_{y \in 2A_\ell^k} \int_{x\in Q_k} |\nabla_x G(x,y)|^2 \,dx \,dy \right)^p. \no
\er
Hence,
\br
&& \sum_{j  \in A^k_+} \Big \langle \int_{Q_j} |\nabla \tilde w_k|^{2p} \,dx  \Big \rangle^{\frac{1}{2p}} \no \\
&& \quad \quad \quad \leq  C \sum_{\ell=1}^{N} |A_\ell^k|^{(2p-1)/(2p)} \left( \Big \langle \int_{A_\ell^k} |\nabla \tilde w_k|^{2p} \,dx \Big \rangle \right)^{1/(2p)} \no \\
&& \quad \quad \quad  \leq  C\sum_{\ell=1}^{N} (2^{\ell})^{d(2p-1)/(2p)} (2^\ell)^{(d -p(2 + d))/(2p) } \Big \langle\int_{y \in 2A_\ell^k} \int_{x\in Q_k} |\nabla_x G(x,y)|^2 \,dx \,dy \Big \rangle^{1/2} \no \\
&& \quad \quad \quad  =  C  \sum_{\ell=1}^{N} (2^{\ell})^{d/2 - 1} \Big \langle \int_{y \in 2A_\ell^k} \int_{x\in Q_k} |\nabla_x G(x,y)|^2 \,dx \,dy \Big \rangle^{1/2}. \label{diadic1}
\er
By stationarity we have 
\[
\Big \langle \int_{y \in Q_j} \int_{x \in Q_k} |\nabla_x G(x,y)|^2 \,dx \,dy \Big \rangle = \Big \langle \int_{y \in Q_0} \int_{x \in Q_{k-j}} |\nabla_x G(x,y)|^2 \,dx \,dy \Big \rangle.
\]
Therefore,
\br
\Big \langle \int_{y \in 2A_\ell^k} \int_{x\in Q_k} |\nabla_x G(x,y)|^2 \,dx \,dy \Big \rangle & = & \Big \langle \int_{y \in Q_0} \int_{x\in 2A_\ell^0} |\nabla_x G(x,y)|^2 \,dx \,dy \Big \rangle. \label{gradGrearrange}
\er
The point here is that the integral in $x$ is now over the annulus $A^0_\ell$ of diameter $O(2^\ell)$, rather than over the unit cube.

For $d \geq 3$, we combine (\ref{gradGrearrange}) with Cacciopoli's inequality to $x \mapsto G(x,y)$. The result is:
\br
\Big \langle \int_{y \in 2A_\ell^k} \int_{x\in Q_k} |\nabla_x G(x,y)|^2 \,dx \,dy \Big \rangle 
& \leq & C (2^{\ell})^{-2} \Big \langle \int_{y \in Q_0} \int_{x\in 2A_\ell^0} |G(x,y)|^2 \,dx \,dy \Big \rangle \no \\
&& + C |D_L|^{-1} \Big \langle \int_{y \in Q_0} \int_{x\in 2A_\ell^0} |G(x,y)| \,dx \,dy \Big \rangle. \no
\er
By Lemma \ref{lem:udecay}, we have a uniform decay estimates for $|G(x,y)| \leq C dist(x,y)^{2-d}$ for $d \geq 3$. Therefore,
\[
\Big \langle \int_{y \in 2A_\ell^k} \int_{x\in Q_k} |\nabla_x G(x,y)|^2 \,dx \,dy \Big \rangle  \leq  C (2^{\ell})^{-2} (2^{\ell})^{d + 2(2 - d)}  +   CL^{-d} (2^\ell)^d (2^\ell)^{2 - d} \leq (2^\ell)^{2 - d}. 
\]
So, returning to (\ref{diadic1}), we obtain
\br
\sum_{j  \in A^k_+} \Big \langle \int_{Q_j} |\nabla \tilde w_k|^{2p} \,dx  \Big \rangle^{1/2p}  \leq  C  \sum_{\ell=1}^{N} (2^\ell)^{d/2 - 1 } (2^\ell)^{1 - d/2} = O(\log L). \label{tildewsum}
\er

In the case $d = 2$, we apply Lemma \ref{lem:ugrad2d} directly to (\ref{gradGrearrange}) and conclude
\[
\Big \langle \int_{y \in 2A_\ell^k} \int_{x\in Q_k} |\nabla_x G(x,y)|^2 \,dx \,dy \Big \rangle  \leq  C. 
\]
So, returning to (\ref{diadic1}), we still obtain
\br
\sum_{j  \in A^k_+} \Big \langle \int_{Q_j} |\nabla \tilde w_k|^2 \,dx \Big \rangle^{1/2}  \leq  C   \sum_{\ell=1}^{O(\log L)} (2^\ell)^{d/2 - 1 }  =  O(\log L). \label{tildewsum2d}
\er
This completes the proof of Lemma \ref{lem:wksum}. \hfill \qed

\vspace{0.2in}

Now we combine (\ref{minkow1}), (\ref{minkow2}), (\ref{genHolderapp0}), (\ref{genHolderapp}), (\ref{momentequiv}) and Lemma \ref{lem:wksum} to conclude that
\br
&& \sum_k \expE  \left[  \left |  \sum_{j \in D_L}  (\Delta_k \Delta_j \Gamma(Z)) \overline{\Delta_j \Gamma(Z^A)} \right|^2 \right] \no \\
&& \quad \quad \quad \quad \leq \sum_{k} \left(   \sum_{j \in A^k_0 \cup A^k_1}  \langle |\Delta_k \Delta_j \Gamma(Z) \overline{\Delta_j \Gamma(Z^A)}|^2 \rangle^{1/2} +  \sum_{\ell=2}^N\sum_{j \in A^k_\ell}  \langle |\Delta_k \Delta_j \Gamma(Z) \overline{\Delta_j \Gamma(Z^A)}|^2 \rangle^{1/2} \right)^2 \no \\
&& \quad \quad \quad \quad \leq \sum_{k} \left( C L^{-2d} \langle \Phi_0^8\rangle^{1/2} +  C \langle \Phi_0^{8q}\rangle^{1/2q}  L^{-2d} \log L \right)^2 \no \\
&& \quad \quad \quad \quad \leq C L^{-3d}\langle \Phi_0^{8q}\rangle^{1/q}(\log L)^2   \label{twosums2}
\er
holds for $d \geq 2$, for all $L \geq 2$, $\beta \geq 0$. 

Finally, we estimate (\ref{secondsum}).   By Minkowsi's inequality we have
\br
\expE  \left[  \left |  \sum_{j \in D_L}  ( \Delta_j \Gamma(Z))^k \Delta_k  \overline{\Delta_j \Gamma(Z^A)} \right|^2 \right] & \leq & \left( \sum_{j \in D_L} \langle |(\Delta_j \Gamma(Z))^k \Delta_k  \overline{\Delta_j \Gamma(Z^A)}|^2 \rangle^{1/2} \right)^2. \label{minkow3}
\er
Recall the notation (\ref{Aavg}) for the average with respect to sets $A$ not containing index $j$. In particular, the weights $K_{n,A}$ define a probability distribution over the index sets $A$ not containing $j$. By applying Jensen' inequality to (\ref{minkow3}) we obtain
\br
 \left( \sum_{j \in D_L} \langle |(\Delta_j \Gamma(Z))^k \Delta_k  \overline{\Delta_j \Gamma(Z^A)}|^2 \rangle^{1/2} \right)^2 & = &  \left( \sum_{j \in D_L} \langle |   \overline{ (\Delta_j \Gamma(Z))^k(\Delta_k (\Delta_j \Gamma))(Z^A,Z',Z'')	}|^2 \rangle^{1/2} \right)^2  \no \\
 & \leq &  \left( \sum_{j \in D_L} \langle \overline{ |(\Delta_j \Gamma(Z))^k (\Delta_k (\Delta_j \Gamma))(Z^A,Z',Z'')|^2} \rangle^{1/2} \right)^2 \no \\
& = &  \left( \sum_{j \in D_L} \langle \! \langle  |(\Delta_j \Gamma)^k(Z,Z',Z'') (\Delta_k (\Delta_j \Gamma))(Z^A,Z',Z'')|^2 \rangle\! \rangle^{1/2} \right)^2, \no 
\er
where we have introduced the notation
\br
\davg{H_A(Z,Z',Z'')} & = & \expE[\overline{H_A(Z,Z',Z'')}]  \no \\
& = & \expE[\sum_{\substack{A \subset [n] \\ j \notin A}} K_{n,A} H_A(Z,Z',Z'')]  =  \sum_{\substack{A \subset [n] \\ j \notin A}}  K_{n,A}  \expE[H_A(Z,Z',Z'')].
\er
The rest proceeds exactly as in the proof of (\ref{twosums2})), the only difference being the following modification of Lemma \ref{lem:genHolder}:

\begin{lem} 
For $p > 1$ there is a constant $C_p$ such that if $k,j \in D_L$ with $|k - j| > 2\tau$ then
\br
&& L^{4d} \davg{ |(\Delta_k \Delta_j \Gamma)(Z^A,Z',Z'') (\Delta_j \Gamma)^k(Z,Z',Z'')|^2}  \leq   C_p \langle \Phi_0^{4q} \rangle^{\frac{4}{2q}} \left\langle (\int_{B_\tau(j)} |\nabla \tilde w_k|^{2} \,dx)^p \right \rangle^{1/p}  \label{matvecjkA}
\er
where $q = 2p/(p-1)$. Also, there is a constant $C$ such that 
\be
L^{4d}\davg{ |(\Delta_k \Delta_j \Gamma)(Z^A,Z',Z'') (\Delta_j \Gamma)^k(Z,Z',Z'')|^2} \leq  C \langle\Phi_0^8\rangle \label{anykjfboundA}
\ee
holds for all $j,k \in D_L$, $L > 1$, $\beta \geq 0$. 
\end{lem}

\vspace{0.2in}

{\bf Proof:} The proof is almost identical to that of Lemma \ref{lem:genHolder}.  We only need to observe that, for any pair of indices $j,k \in D_L$ and any set $A \subset D_L \cap \mathbb{Z}^d$, if $g(Z,Z',Z'')$ denotes any of the random variables $\hat \Phi_j^2(Z)$, $\hat \Phi_j^2(Z^j)$, $\hat \Phi_j^2(Z^{k})$, $\hat \Phi_j^2(Z^{jk})$, $\hat \Phi_k^2(Z)$, $\hat \Phi_k^2(Z^{j})$, $\hat \Phi_k^2(Z^k)$, or $\hat \Phi_k^2(Z^{jk})$, then $g(Z,Z',Z'')$ and $g(Z^A,Z',Z'')$ have the same distribution. In particular, for any power $p$,
\be
\davg{g(Z^A,Z',Z'')^p} = \davg{g(Z,Z',Z'')^p} = \avg{g(Z,Z',Z'')^p} = \avg{\hat \Phi_0^{2p}}.  \label{gavgequiv}
\ee
Similarly, the random variables $(\Delta_j \Gamma)^k(Z,Z',Z'')$ and $\Delta_j \Gamma(Z,Z')$ have the same distribution. As before, by Lemma \ref{lem:DeljDelkG} and Lemma \ref{lem:DeljG}, we have
\br
L^{2d} |(\Delta_k \Delta_j \Gamma)(Z,Z',Z'')|^2 & \leq & C \left(\hat \Phi_j^2(Z) + \hat \Phi_j^2(Z^j)+ \hat \Phi_j^2(Z^{k}) + \hat \Phi_j^2(Z^{jk})\right) \no \\
&& \quad \quad \times \left(\hat \Phi_k^2(Z) + \hat \Phi_k^2(Z^j)+ \hat \Phi_k^2(Z^{k}) + \hat \Phi_k^2(Z^{jk})\right)  \no
\er
and
\[
L^{2d} \langle  | (\Delta_j \Gamma)^k(Z,Z',Z'')|^{2q} \rangle^{1/q} = L^{2d} \langle  | \Delta_j \Gamma(Z,Z')|^{2q} \rangle^{1/q}    \leq C \langle \hat \Phi_0^{4q} \rangle^{\frac{1}{q}}. 
\]
Therefore, as in the proof of Lemma \ref{lem:genHolder}, by applying the generalized H\"older inequality and (\ref{gavgequiv}) we obtain
\br
L^{2d} \davg{ |(\Delta_k \Delta_j \Gamma)(Z^A,Z',Z'') (\Delta_j \Gamma)^k(Z,Z',Z'')|^2 }   &\leq & C \avg{ \hat \Phi_0^{8} }^{\frac{1}{4}} \avg{\hat \Phi_0^{8} }^{\frac{1}{4}}  \langle  |(\Delta_j \Gamma)^k|^{4} \rangle^{1/2} \no \\
&  \leq & C \langle \hat \Phi_0^{8} \rangle^{\frac{1}{4}} \langle \hat \Phi_0^{8} \rangle^{\frac{1}{4}}  L^{-2d} \langle \hat \Phi_0^{8} \rangle^{\frac{1}{2}}.  \label{genHolderboundA}
\er
This and (\ref{momentequiv}) imply (\ref{anykjfboundA}).

If $dist(k,j) > 2 \tau$, Lemma \ref{lem:DeljDelkG} tells us that
\br
L^{2d} |\Delta_k \Delta_j \Gamma(Z)|^2 & \leq & C \left(\hat \Phi_j^2(Z) + \hat \Phi_j^2(Z^j)+ \hat \Phi_j^2(Z^{k}) + \hat \Phi_j^2(Z^{jk})\right) \no \\
&& \quad \quad \times \left(  \hat \Phi_k^2(Z^k) \int_{B_\tau(j)} |\nabla \tilde w_k|^2 \,dx +  \Phi_k^2(Z^{jk}) \int_{B_\tau(j)} |\nabla \tilde w_k^j|^2\,dx \right), \label{DeljDelkG2wA}
\er
where
\[
\tilde w_k = \hat \Phi_k(Z^k)^{-1} \Delta_k \phi, \quad \quad \tilde w_k^j = \hat \Phi_k(Z^{jk})^{-1} \Delta_k \phi^j.
\]
Let $p > 1$, let $q = 2p/(p-1)$ so that $\frac{1}{p}  + \frac{1}{2q} + \frac{1}{2q} + \frac{1}{q}= 1$. Then by (\ref{DeljDelkG2wA}) and the generalized H\"older inquality,
\br
&& L^{2d} \davg{ |(\Delta_k \Delta_j \Gamma)(Z^A,Z',Z'') (\Delta_j \Gamma)^k(Z,Z',Z'') |^2 }   \no \\
& & \quad \quad \quad \quad \leq   C \langle \hat \Phi_0^{4q} \rangle^{\frac{1}{2q}} \langle \hat \Phi_0^{4q} \rangle^{\frac{1}{2q}} \davg{g_{jk}(Z^A,Z_k'')^p}^{1/p} \langle  | (\Delta_j \Gamma)^k(Z,Z',Z'') |^{2q} \rangle^{1/q} \no \\
& & \quad \quad \quad \quad \quad +  C \langle \hat \Phi_0^{4q} \rangle^{\frac{1}{2q}} \langle \hat \Phi_0^{4q} \rangle^{\frac{1}{2q}} \davg{ g_{jk}(Z^{j\cup A},Z_k'')^p }^{1/p} \langle  | (\Delta_j \Gamma)^k(Z,Z',Z'') |^{2q} \rangle^{1/q},  \label{genHolderboundA2}
\er
where
\[
g_{jk}(Z,Z_k'') = \int_{B_\tau(j)} |\nabla_x \tilde w_k(x,Z,Z_k'')|^{2} \,dx.
\]
On the other hand, $g_{jk}(Z,Z_k'')$ and $g_{jk}(Z^A,Z_k'')$ and $g_{jk}(Z^{j \cup A},Z_k'')$ all have the same distribution. Hence
\[
\davg{ g_{jk}(Z^{j\cup A},Z_k'')^p }^{1/p}  = \davg{ g_{jk}(Z^{A},Z_k'')^p }^{1/p}  =\Big \langle (\int_{B_\tau(j)} |\nabla \tilde w_k|^{2} \,dx)^p \Big \rangle^{1/p}.
\]
We conclude that
\br
&& L^{2d} \davg{ |(\Delta_k \Delta_j \Gamma)(Z^A,Z',Z'') (\Delta_j \Gamma)^k(Z,Z',Z'') |^2 }   \no \\
& & \quad \quad \quad \quad \leq   C \langle \hat \Phi_0^{4q} \rangle^{\frac{1}{2q}} \langle \hat \Phi_0^{4q} \rangle^{\frac{1}{2q}} \Big \langle (\int_{B_\tau(j)} |\nabla \tilde w_k|^{2} \,dx)^p \Big \rangle^{1/p} \langle  |\Delta_j \Gamma(Z)|^{2q} \rangle^{1/q}  \no \\
& & \quad \quad \quad \quad \leq   C \langle \hat \Phi_0^{4q} \rangle^{\frac{1}{2q}} \langle \hat \Phi_0^{4q} \rangle^{\frac{1}{2q}} \Big \langle (\int_{B_\tau(j)} |\nabla \tilde w_k|^{2} \,dx)^p \Big \rangle^{1/p} L^{-2d} \langle  \Phi_0^{4q} \rangle^{1/q} 
\er
which implies (\ref{matvecjkA}). \hfill \qed

\vspace{0.2in}

With this modification of Lemma \ref{lem:genHolder}, we proceed exactly as in the proof of (\ref{twosums2}) to obtain the bound
\be
\sum_k \expE  \left[  \left |  \sum_{j \in D_L}  ( \Delta_j \Gamma)^k(Z))^k \Delta_k \overline{ \Delta_j \Gamma(Z^A)} \right|^2 \right]  \leq C \langle  \Phi_0^{8q}\rangle^{1/q} L^d L^{-4d} (\log L)^2. \label{twosums3}
\ee
By combining Theorem \ref{theo:normalapprox} with (\ref{sig3sum}), (\ref{twosums2}), and (\ref{twosums3}) we conclude that
\be
d_{\mathcal{W}}\left(\frac{\Gamma_{L,\beta} - m_{L,\beta}}{\sigma_{L,\beta}},Z \right) \leq C \frac{L^{-2d}}{\sigma^3} \expE[  \Phi_0^6] + C \frac{L^{-3d/2} \log(L)}{\sigma^2} \expE[  \Phi_0^{8q}]^{\frac{1}{2q}}. \label{dWZsig}
\ee
for all $d \geq 2$. The exponent $q > 2$ is the H\"older conjugate of $p$, where $2p$ is the exponent from Meyers' estimate.  This concludes the proof of Theorem \ref{theo:cltmain}.

\vspace{0.2in}

\section{Estimates for the periodic Green's function} \label{sec:Greenbound}

\subsection*{$d \geq 3$: Proof of Lemma \ref{lem:udecay}} 

Here we follow ideas used to prove a uniform decay estimate for Green's functions in $\Rm^d$, as in Theorem 1.1 of \cite{GW} and Lemma 2.8 of \cite{GO}; the difference here is the periodicity, so we include a proof for completeness. Let $y \in D_L$ and let $u(x) = G(x,y)$ be the periodic Green's function, which satisfies
\be
- \nabla \cdot (a \nabla u) + \beta u = \delta_y - |D_L|^{-1}, \quad x \in D_L. \label{ueqn}
\ee
in the weak sense. Suppose that $u$ also satisfies
\be
|\{ x \in D_L\;|\; u(x) > 0 \}| \leq \frac{1}{2} |D_L|. \label{unormcond}
\ee
(If this is not the case, then we could apply the same argument to the function $-u$ instead.) Then, for any $k > 0$, the function $u_k(x) = \max(0,\min(u,k))$ satisfies
\be
|\{ x \in D_L\;|\; u_k(x) \neq 0 \}| = |\{ x \in D_L\;|\; u(x) > 0 \}| \leq \frac{1}{2} |D_L|. \label{ukvanish}
\ee
Since $\norm{u_k}_\infty \leq k$, we observe that $u_k$ satisfies
\[
\int_{D_L} \nabla u_k \cdot a \nabla u_k \,dx = \int_{D_L} \nabla u \cdot a \nabla u_k \,dx \leq   - \beta \int_{D_L} u u_k \,dx + 2k  \leq 2k.
\]
Therefore,
\be
\int_{D_L} |\nabla u_k|^2 \,dx  \leq  2 k/a_* . \label{ukgrad}
\ee

\commentout{
This is a Poincar\'e inequality, albeit nonstandard. Because $q$ is the critical exponent, a simple rescaling shows that (\ref{ukSob}) is equivalent to 
\be
\left( \int_{D_1} |v|^q \,dx \right)^{1/q} \leq C \left( \int_{D_1} |\nabla v|^2 \,dx \right)^{1/2} \label{ukSob2}
\ee
for functions $v \in H^1_{per}(D_1)$ satisfying $|\{ x \in D_1\;|\; v(x) = 0 \}| \geq 1/2$. By the Sobolev imbedding theorem, we know that
\be
\left( \int_{D_1} |v|^q \,dx \right)^{1/q} \leq C \left(\int_{D_1} v^2 \,dx +  \int_{D_1} |\nabla v|^2 \,dx \right)^{1/2} \label{ukSob3}
\ee
must hold for all $v \in H^1_{per}(D_1)$. On the other hand, by the usual Poincar\'e inequality (using the compactness of the embedding of $H^1(D_1)$ into $L^2(D_1)$) we know there is a constant $C$ such that 
\be
\int_{D_1} v^2 \,dx \leq C \int_{D_1} |\nabla v|^2 \,dx \label{ukSob4}
\ee
holds for all functions $v \in H^1_{per}(D_1)$ which also satisfy $|\{ x \in D_1\;|\; v(x) = 0 \}| \geq 1/2$ (for example, see Lemma 4.8 of \cite{HL}). So, (\ref{ukSob2}) follows from (\ref{ukSob3}) and (\ref{ukSob4}). 
}

Considering (\ref{ukvanish}), we know there is a constant $C$, independent of $k$, $L$ and $\beta$, such that
\be
\left( \int_{D_L} |u_k|^q \,dx \right)^{1/q} \leq C \left( \int_{D_L} |\nabla u_k|^2 \,dx \right)^{1/2} \label{ukSob}
\ee
where $q = 2d/(d-2)$ is the critical Sobolev exponent. By scaling, this is a consequence of the Sobolev imbedding theorem and the Poincar\'e inequality for functions $v \in H^1_{per}(D_1)$ which also satisfy $|\{ x \in D_1\;|\; v(x) = 0 \}| \geq 1/2$ (for example, see Lemma 4.8 of \cite{HL}). By applying Chebychev's inequality, then (\ref{ukSob}) and (\ref{ukgrad}), we obtain the estimate
\br
| \{ x \in D_L \;|\; u(x) \geq k \}|  =  | \{ x \in D_L \;|\; u_k(x) \geq k \}| & \leq & k^{-q}\int_{D_L} |u_k|^q \,dx\leq  C k^{-q/2}. \label{uweakq}
\er
This is a weak-$L^{p}(D_L)$ estimate on $u^+ = \max(u,0)$, for $p = q/2 = d/(d-2)$:
\be
\norm{u^+}_{L^{p}_W(D_L)} = \sup_{t > 0} \;\;t \;| \{ x \in D_L \;|\; |u^+(x)| > t \}|^{1/p} \leq C,\label{uweakp}
\ee  
where the constant $C$ is independent of $L$ and $\beta \geq 0$.  

Now let $\alpha \in (1, p)$, $x_0 \in D_L$, $R < dist(x_0,y)$. The weak bound (\ref{uweakp}) implies that $u^+ \in L^\alpha(B_R(x_0))$. By using the identity
\[
\int_{B_R} |u^+|^{\alpha}\,dx = \alpha \int_0^\infty t^{\alpha-1} |\{ x \in B_R \;|\; u^+(x) \geq t \}| \,dt \leq |B_R| s^\alpha +  \alpha \int_s^\infty t^{\alpha-1} |\{ x \in B_R \;|\; u^+(x) \geq t \}| \,dt 
\] 
and optimizing in $s$, we see that 
\be
\norm{u^+}_{L^{\alpha}(B_R)} \leq C \left( \frac{p}{p-\alpha}\right)^{1/\alpha} |B_R|^{\frac{p-\alpha}{p\alpha}}, \label{uweakstrong}
\ee
where the constant $C$ depends on $\alpha$ and $p$, but not on $L$ or $R$ or $\beta \geq 0$. Since $-\nabla \cdot (a \nabla u) + \beta u = -|D_L|^{-1} $ in $B_R$, the estimates of De Giorgi and Moser give us a bound on $u^+(x)$ in terms of $\norm{u^+}_{L^\alpha(B_R(x_0))}$. Specifically, Theorem 4.1 of \cite{HL} (or Theorem 8.17 of \cite{GT}) implies that $u$ is locally bounded and satisfies:
\br
\sup_{x \in B_{R/2}(x_0)} u^+(x) & \leq & C  R^{-d/\alpha} \left( \int_{B_R} (u^+(y))^{\alpha} \,dy \right)^{1/\alpha} + C R^{2}|D_L|^{-1}, \label{ulocbound}
\er
with a constant $C$ that depends only on $d$, $a_*$, $a^*$, and $\alpha$. Note that in Theorem 4.1 of \cite{HL}, the constant depends on $|\beta|R^2$. However, it is easy to see from the proof (method 1) that if $\beta$ is known to be non-negative, then the bound is independent of $\beta$, so the same bound holds under rescaling (as in Theorem 4.14 of \cite{HL}).

By combining (\ref{uweakstrong}) and (\ref{ulocbound}) we have
\br
\sup_{x \in B_{R/2}(x_0)} u^+(x)  \leq  C  R^{-d/\alpha} |B_R|^{\frac{p-\alpha}{p\alpha}}   +  C R^{2}|D_L|^{-1}  \leq  C R^{-d/p}  + C R^{2}L^{-d} \leq C R^{2-d}, \no
\er
where the constant $C$ depends on the dimension, but not on $L$, $\beta \geq 0$, $R$. In particular,
\be
u^+(x) \leq C \left( dist(x,y)\right)^{2-d}. \label{upupper}
\ee

Now, assuming (\ref{unormcond}) holds for $u$ (otherwise, replace $u$ by $-u$), let us choose $r \leq 0$ such that both
\[
|\{ x \in D_L\;|\; u(x) > r \}| \leq \frac{1}{2}|D_L|
\]
and
\[
|\{ x \in D_L\;|\; u(x) < r \}| \leq \frac{1}{2}|D_L|
\]
hold. Consider the function $\bar u = r - u$ which satisfies
\[
- \nabla \cdot (a \nabla \bar u) +\beta \bar u = - \delta_y + |D_L|^{-1} - \beta |r|
\]
and 
\[
|\{ x \in D_L\;|\; \bar u(x) > 0 \}| \leq \frac{1}{2}|D_L|.
\]
To the functions $\bar u_k = \max(0,\min(\bar u,k))$ and $\bar u^+ = \max(0,\bar u)$ we apply the same argument used to obtain (\ref{upupper}). The result is:
\be
\bar u^+(x) \leq C \left( dist(x,y)\right)^{2-d}. \label{ubarupper}
\ee
In deriving (\ref{ulocbound}) for $\bar u^+$, we must use the fact that $\bar u$ is a subsolution of $- \nabla \cdot (a \nabla \bar u) +\beta \bar u = |D_L|^{-1}$ away from $y$, since $-\beta |r| \leq 0$. That is,
\[
\int_{B_R}  \nabla \varphi  \cdot a \nabla \bar u+ \beta \bar u \varphi \,dx \leq |D_L|^{-1} \int_{B_R} \varphi \,dx
\]
holds for all $\varphi \in H^1_{0}(B_R)$ which satisfy $\varphi \geq 0$. Thus, Theorem 4.1 of \cite{HL} (or Theorem 8.17 of \cite{GT}) still applies. Apart from this detail, the argument is identical.  By combining (\ref{upupper}) and (\ref{ubarupper}) we obtain
\be
r - C \left( dist(x,y)\right)^{2-d} \leq r - \bar u^+(x) \leq u(x) \leq u^+(x) \leq C \left( dist(x,y)\right)^{2-d} \label{rsqueeze}
\ee

On the other hand, (\ref{upupper}) implies that
\[
\int_{D_L} u^+(x) \,dx \leq C L^2.
\]
We combine this with the fact that $\int_{D_L} u \,dx = 0$ to conclude that
\br
0 & = & \int_{D_L} u^+(x) \,dx + \int_{\{ r < u \leq 0\}} u(x) \,dx + \int_{\{ u \leq r \}} u(x) \,dx \no \\
& \leq&  CL^2 + r |\{ x \in D_L \;|\; u(x) \leq r \}| \no \\
& = &   CL^2 + r |\{ x \in D_L \;|\; \bar u(x) \geq 0 \}|  \leq C L^2 + r L^{d}/2. \no
\er
Hence $|r| \leq 2CL^{2 - d}$. Combining this with (\ref{rsqueeze}) we obtain $|u(x)| \leq C dist(x,y)^{2-d}$, as desired. \hfill \qed

\subsection*{$d=2$: Proof of Lemma \ref{lem:ugrad2d}}

Lemma \ref{lem:ugrad2d} relies on the following oscillation estimate, which is a version of Lemma 2.8(i) of \cite{GO}:

\begin{lem} \label{lem:udecay2d} 
Let $d =2$. For any $q \geq 1$, there is a constant $C > 0$ such that
\[
R^{-2} \int_{B_R(x_0)} |G(x,y)- \bar G_R(y)|^q \,dx \leq C
\]
holds for all $x_0 \in D_L$, $y \in D_L \setminus B_{2R}(x_0)$, $R > 0$, $L > 1$ and $\beta \geq 0$, where $\bar G_R(y)$ is the average of $G(\cdot,y)$ over the ball $B_R(x_0)$.
\end{lem}

{\bf Proof of Lemma \ref{lem:udecay2d}:} This is proved as in Lemma 2.8 of \cite{GO} (see part (i), Step 2) for the free-space Green's function (see Step 2 in the proof therein); here we include the proof for completeness. Fix $y \in D_L$. Let $u(x) = G(x,y)$, which satisfies
\[
- \nabla \cdot (a \nabla u ) + \beta u = \delta_y - |D_L|^{-1}.
\]
Let $\overline{u_R}$ be the average of $u$ over the ball $B_R$. Without loss of generality, suppose $\overline{u_R} \geq 0$. For $k \geq 0$, define
\[
u_k = \max( \min(u, \overline{u_R} + k), \overline{u_R} - k).
\]
We claim that
\be
\int_{D_L} |\nabla u_k|^2  \,dx \leq \frac{2k}{a_*}. \label{gradukk}
\ee
To see this, observe that for any constant $c \in \Rm$, 
\br
\int_{D_L} \nabla u_k \cdot a \nabla u_k  \,dx & = & \int_{D_L} \nabla (u_k + c) \cdot a \nabla u  \,dx  \no \\
& = & u_k(y) - |D_L|^{-1}\int_{D_L} u_k(x) \,dx  - \beta \int_{D_L} u  (u_k + c) \,dx. \label{ukeqn}
\er
If $\overline{u_R} \in [0,k]$, let $c = 0$. Then $u(x)(u_k(x) + c) \geq 0$ at every point $x \in D_L$.  Hence
\be
\beta \int_{D_L} u  (u_k + c) \,dx \geq 0. \label{ukucint}
\ee
Therefore, (\ref{gradukk}) follows from (\ref{ukeqn}). If $\overline{u_R} > k$, let $c = k - \overline{u_R}$. Then $u_k + c \geq 0$. Also, $u(x) > \overline{u_R} - k > 0$ must hold wherever $(u_k(x) + c) > 0$.  Hence (\ref{ukucint}) still holds. Moreover, $0 \leq u_k(x) + c \leq 2k$, so again (\ref{gradukk}) follows from (\ref{ukeqn}).

Now let $v(x) = u(x) - \overline{u_R}$. Let $v_k(x) = \max(\min(v(x),k),-k) = u_k(x) - \overline{u_R}$. Let $\overline{v_R}$ and $\overline{v_{k,R}}$ be the average of $v$ and $v_k$ over $B_R$, respectively. Hence $\overline{v_R} = 0$. Then the goal is to bound
\br
\left( R^{-2} \int_{B_R} |v|^q \,dx \right)^{1/q} & = &  \left( R^{-2} \int_{B_R \cap \{ |v| \leq k\}} |v_k|^q \,dx + R^{-2} \int_{B_R \cap \{ |v| > k\}} |v|^q \,dx  \right)^{1/q} \no \\
& \leq & C \left( R^{-2} \int_{B_R \cap \{ |v| \leq k\}} |v_k - \overline{v_{k,R}}|^q \,dx\right)^{1/q} + C |\overline{v_{k,R}}| \no \\
& & + C \left(  R^{-2} \int_{B_R \cap \{ |v| > k\}} |v|^q \,dx  \right)^{1/q}. \no
\er
Since $\overline{v_R} = 0$, we have
\[
|\overline{v_{k,R}}| \leq 2 \left( R^{-2} \int_{B_R \cap \{|v| \geq k\}} |v|^q \,dx \right)^{1/q}.
\]
Therefore,
\br
\left( R^{-2} \int_{B_R} |v|^q \,dx \right)^{1/q} & \leq & C \left( R^{-2} \int_{B_R} |v_k - \overline{v_{k,R}}|^q \,dx\right)^{1/q} \no \\
& & + C \left(  R^{-2} \int_{B_R \cap \{ |v| > k\}} |v|^q \,dx  \right)^{1/q}. \label{Iqfirst}
\er
By the Sobolev inequality and then (\ref{gradukk}), we know that for any $s \in [1, \infty)$ there is a constant $C_s$ (depending only on $s$) such that
\[
\left( R^{-2} \int_{B_R} |v_k - \overline{v_{k,R}}|^s \,dx \right)^{1/s}  \leq  C_s \left( \int_{B_R} |\nabla v_k|^2 \,dx \right)^{1/2} =  C_s \left( \int_{B_R} |\nabla u_k|^2 \,dx \right)^{1/2} \leq C k^{1/2}. 
\]
To estimate the last integral appearing in (\ref{Iqfirst}) we use
\[
\int_{B_R \cap \{ |v| > k\}} |v|^q \,dx  =  \int_0^\infty q t^{q - 1}|\{ |v| \geq \max(t,k) \}| \,dt \leq  |\{ |v| \geq k \}| k^{q} + q \int_k^\infty t^{q - 1} |\{ |v| \geq t \}| \,dt,
\]
and
\[
|\{ |v| \geq k \}| \leq  |\{ |v_k| \geq k \}| \leq k^{-s}  \int_{B_R} |v_k|^s \,dx.
\]
Let $s > 2q$. Then
\br
\int_{B_R} |v_k|^s \,dx & \leq & C  \int_{B_R} |v_k - \overline{v_{k,R}}|^s \,dx  +  C R^2 (\overline{v_{k,R}})^s \no \\
& \leq & C  \int_{B_R} |v_k - \overline{v_{k,R}}|^s \,dx  +  C R^2\left(  R^{-2} \int_{B_R} |v|^q \,dx \right)^{s/q} \no \\
& \leq & C  R^2 k^{s/2}   +  C R^2\left(  R^{-2} \int_{B_R} |v|^q \,dx \right)^{s/q}. \no
\er
So, if $I_q = (R^{-2} \int_{B_R} |v|^q \,dx)^{1/q}$, we have
\[
\int_{B_R} |v_k|^s \,dx  \leq   C  R^2 k^{s/2}   +  C R^2 I_q^s
\]
and $|\{ |v| \geq k \}| \leq   R^2 k^{-s/2}   +  C k^{-s} R^2  I_q^s$.

Combining these bounds and returning to (\ref{Iqfirst}), we obtain
\br
I_q &  \leq & C k^{1/2} + C R^{-2/q} \left(  |\{ |v| \geq k \}| k^{q } + q \int_k^\infty t^{q - 1} |\{ |v| \geq t \}| \,dt \right)^{1/q} \no \\
& \leq & C k^{1/2} + C R^{-2/q} \left(  R^2 k^{-s/2}k^q   +  C k^{-s} R^2  I_q^s k^{q}\right)^{1/q} \no \\
& & + C R^{-2/q} \left(q \int_k^\infty t^{q - 1} |\{ |v| \geq t \}| \,dt \right)^{1/q} \no \\
& \leq & C k^{1/2} + C R^{-2/q} \left(  R^2 k^{-s/2}k^q   +  C k^{-s} R^2  I_q^s k^{q}\right)^{1/q} \no \\
& & + C R^{-2/q} \left(q \int_k^\infty t^{q - 1} (R^2 t^{-s/2}   +  C t^{-s} R^2  I_q^s) \,dt \right)^{1/q} \no \\
& \leq & C k^{1/2} + C k^{1 - s/(2q)} + C I_q^{s/q} k^{1 - s/q}.
\er
By choosing $k = \alpha I_q$ with $\alpha > 0$ sufficiently large, we see that this implies $I_q \leq C$. \hfill \qed

\vspace{0.2in}

Now we continue with the proof of Lemma \ref{lem:ugrad2d}. By assumption, $dist(x_0,y) > 2R$. Let $\varphi$ be a smooth function supported in $B_{2R}(x_0)$ and satisfying: $0 \leq \varphi(x) \leq 1$ for all $x$, $\varphi(x) = 1$ for $x \in B_{R}(x_0)$, and $|\nabla \varphi| \leq C/R$. Applying Lemma \ref{lem:CaccioGen} to $u(x) = G(x,y)$ with this choice of $\varphi$, we conclude
\br
\int_{B_{R}} \abs{\nabla u}^2  \, dx  & \leq &   K_1 |D_L|^{-1} \int_{B_{2R}} (u - b) \varphi^2 \,dx   -K_1  \beta \int_{B_{2R}} u(u-b) \varphi^2 \,dx \no \\
& & + K_2  \int_{B_{2R}} \abs{\nabla \varphi}^2 (u- b)^2\,dx.  
\er
If we choose 
\[
b = \left(\int_{B_{2R}} \varphi^2 \,dx \right)^{-1} \int_{B_{2R}} u \varphi^2 \,dx,
\]
then Jensen's inequality implies
\[
 \int_{B_{2R}} u(u - b) \varphi^2 \,dx \geq 0.
\]
Therefore, since $\beta \geq 0$,
\br
\int_{B_{R}} \abs{\nabla u}^2  \, dx  & \leq &  K_1 |D_L|^{-1} \int_{B_{2R}} (u - b) \varphi^2 \,dx +     K_2  \int_{B_{2R}} \abs{\nabla \varphi}^2 (u- b)^2\,dx \no \\
& &  \leq K_1 R^{-2} \int_{B_{2R}} |u - b| \,dx + C R^{-2} \int_{B_{2R}} (u- b)^2\,dx. \label{caccioGenbound2}
\er
On the other hand, if $\bar u$ denotes the average of $u(x)$ over $B_{2R}(x_0)$, we know from Lemma \ref{lem:udecay2d} that
\[
R^{-2} \int_{B_{2R}} |u - \bar u| \,dx \leq C.
\]
Hence
\[
|\bar u - b| \leq  \left(\int_{B_{2R}} \varphi^2 \,dx \right)^{-1} \int_{B_{2R}} |\bar u - u(x)| \varphi^2 \,dx \leq C R^{-2}  \int_{B_{2R}} |\bar u - u(x)| \varphi^2 \,dx \leq C.
\]
Applying Lemma \ref{lem:udecay2d} again, we obtain
\[
R^{-2} \int_{B_{2R}} (u- b)^2\,dx \leq C R^{-2} \int_{B_{2R}} (u- \bar u)^2\,dx + C R^{-2} \int_{B_{2R}} (\bar u- b)^2\,dx \leq C.
\]
Similarly,
\[
R^{-2} \int_{B_{2R}} |u- b| \,dx \leq C R^{-2} \int_{B_{2R}} |u- \bar u| \,dx + C R^{-2} \int_{B_{2R}} (\bar u- b)\,dx \leq C.
\]
In view of (\ref{caccioGenbound2}) and the fact that $C$ is independent of $R$, $L$ and $\beta \geq 0$, we have proved the desired result. \hfill \qed

\bibliographystyle{plain}

\bibliographystyle{plain}

\end{document}